\documentclass[a4paper, 12pt]{article}
\usepackage{enumerate, theorem}
\usepackage{amsmath, amsfonts, amssymb}
\usepackage[height=22.5cm, width=15cm]{geometry}
\usepackage[all]{xy}
\usepackage{mathrsfs, yfonts}
\usepackage[utf8]{inputenc}

\newcommand{\Abc}{\tilde{\mathcal{A}}_{conv.}}

\newcommand{\parag}[1]{\paragraph{\sc{#1.}}}

\newtheorem{thm}{Theorem}[subsection]
\newtheorem{defn}[thm]{Definition}
\newtheorem{cor}[thm]{Corollary}
\newtheorem{prop}[thm]{Proposition}
\newtheorem{lemma}[thm]{Lemma}

\begin{document}
\title{Generalized Brieskorn Modules III\\
The algebra $\Abc$}

\author{Daniel Barlet\footnote{Barlet Daniel, Institut Elie Cartan UMR 7502  \newline
Universit\'e de Lorraine, CNRS,  et  Institut Universitaire de France, \newline
BP 239 - F - 54506 Vandoeuvre-l\`es-Nancy Cedex.France. \newline
e-mail : daniel.barlet@univ-lorraine.fr}.}

\maketitle

 {\it A mon ami Alan Huckleberry qui vient de nous quitter, 
ce m\'elange d'Alg\`ebre,  d'Analyse  et de G\'eom\'etrie qui lui aurait peut-être plut}.

\parag{Abstract} In this paper we introduce and study the  "convergent" algebra (containing "a" and "b" and acting on holomorphic germs in "a") which naturally acts on the "generalized Brieskorn modules" associated to the Gauss-Manin connections of  the germs at each point  of the singular set of a holomorphic function on a complex manifold. We generalize to this convergent setting the results previously obtained (see {8], [9],  [15] and [16])  in the formal case, and we  show that, in suitable global situations (for instance when f is projective) we obtain also generalized (geometric) Brieskorn modules. So the question of the relationship between the left module  structure on this algebra (which defines several interesting  filtrations) and the mixte Hodge structure is raised.

\parag{AMS classification}  32 S 25; 32 S 40 ; 34 E 05

\parag{Key words} Asymptotic Expansions ; Period-integral; Generalized Brieskorn Module (GBM);  Convergent (a,b)-Module;  Geometric (convergent) (a,b)-Module.

\tableofcontents

\bigskip

\section*{Introduction}

 Brieskorn modules come from the fundamental paper \cite{[Br.70]} (see also \cite{[S.89]} and \cite{[B.-S. 04]}). The (a,b)-module point of view  appears in singularity theory after the remark, due to Kyoji Saito \cite{[K.S]},  that  the algebra $B$  of Gevrey convergent series in the inverse of the Gauss-Manin connection acts naturally on the Brieskorn module. This defines, in addition to the obvious $A := \mathbb{C}\{s\}$-module structure (here $f = s$), a (non obvious) $B := \mathbb{C}\{\{b\}\}$-module structure, where $b := \int_0^s \simeq \partial_s^{-1}$, with the commutation relation $ab - ba = b^2$ where $a$ is the multiplication by $s$.\\
For quite a long period, I develop the theory of (a,b)-module using the algebra $\widehat{B}[a]$, where $\widehat{B}:= \mathbb{C}[[b]], b := \partial_s^{-1}$,  with the commutation relation
 $$aS(b) = S(b)a + b^2S'(b)$$
 for  $S \in \widehat{B}$, where $S'$ is the usual derivative of $S$ in $\widehat{B}$. See for instance \cite{[B.93]}, \cite{[B.97]}, \cite{[B.06]}, \cite{[B.08]}, \cite{[B.09]} and \cite{[B.22]}.\\
 This point of view highlights the fact that the structure as a B-module is fundamental, placing the A-module structure in second position, contrary to the standard  intuition and usage  in the theory of linear differential systems. \\
 However, while replacing $B$ by $\widehat{B}$ was easily justified by the fact that  the Gauss-Manin  connections  are always regular and the formal completion in “a” implies the formal completion in “b” in the situations at hand, it is nonetheless unsatisfactory  to have only the action of $\mathbb{C}[a]$ and not  an $A$-module structure while such a structure  it is clearly present from the beginning.\\ 
 In my recent papers \cite{[part I]} and \cite{[part II]} I was compelled\footnote{To be able to complete the proofs !} to  use the convergent geometric (a,b)-modules, which are "generalized Brieskorn modules". In these papers, I avoid a part of the difficulty because  it was enough, in order  to complete my proofs,  to work with the algebra $B[a]$. This avoids to face the full real problem which is to  work with the algebra $\Abc$ which contains  both the algebra $B$ and the algebra $A$ (with the standard commutation relation).\\
 The first goal of the present article is to go all the way and   to build the natural framework for the "Generalized Brieskorn Modules" which leads to consider  left $\Abc$-modules which are free finite type $B$-modules. Note that in the "geometric case"  they are also free, finite type $A$-modules which is the standard basic framework for germs of  vector bundles with a connection on a complex curve (see \cite{[D.70]}).\\
 
  So the first aim  of this paper is to introduce the algebra $\Abc$ which is a sub-algebra of the algebra of continuous $\mathbb{C}$-linear endomorphisms of the algebra $A$ which contains both $A$ and $B$ with the suitable actions. Of course the first goal is to show that geometric (convergent) (a,b)-modules are naturally
 left $\Abc$-modules (which are free, finite type $B$-modules) and that, in the geometric case,  this structure is uniquely defined by their left  $B[a]$-module structure, since  $B[a]$-linear  morphisms between geometric generalized Brieskorn modules  extend automatically to $\Abc$-linear morphisms.
  This is obtained by proving a Division Theorem (see Theorem \ref{div.3})  in the algebra $\Abc$.\\
  
    We show also that the algebra $\Abc$ is local, that is to say that a series in (a, b), $u \in \Abc$, with non zero constant term  is invertible in $\Abc$. This property is not true in the algebra $B[a]$ and  allows to show that the action on a geometric convergent  (a, b)-module  of an element in $B[a]$ with a non zero contant term is bijective.\\
   Then, we prove  that $\Abc$ is noetherian (left and right) and is  "almost principal" (see Theorem \ref{24/3}).\\
   
  Our main result  is Theorem \ref{5/10} which is  the generalization to the convergent case  of  Theorem 2.1.1 in \cite{[B.08]}; it  constructs a fonctorial complex of sheaves of 
   $\Abc$-modules on the zero set  $Y$ of a holomorphic function $f$ on a complex manifold $M$, whose cohomology germs at each point $y \in Y$ are, modulo $b$-torsion, the generalized Brieskorn modules at $y$ for $f$. They are convergent  geometric (a,b)-modules (so geometric generalized Brieskorn modules with the terminology introduced in the present article) endowed with their left  $\Abc$-module structure. They represent (up to $b$-torsion)  in the derived category of left  $\Abc$-modules, the Gauss connection of the germs  
   $(f, y)$  for each $y \in Y$, \\
    We conclude this article by proving the finiteness Theorem \ref{finite} which allows to raise the question of the relationship  between the  "global  generalized Brieskorn modules" defined in a suitable global setting, with the corresponding   mixte  Hodge structure on the cohomology of the singular fiber $Y := \{f = 0\}$.

\newpage

\section{The algebra $\Abc$}

\subsection{First properties}

\parag{Notations} In the sequel we shall use the following complex unitary algebras:
\begin{enumerate}
\item The algebras $A_0 := \mathbb{C}[a]$ and $B_0 := \mathbb{C}[b]$. 
\item The algebra $\mathcal{A}_0$ is the algebra of polynomials in $a$ and $b$ inside the algebra $End(\mathbb{C}[s])$ of  linear operators on $A_0$ where $a$ acts as the product by $a$ and $b$ as the primitive without constant. Of course they satisfy $ab - ba = b^2$.
\item The algebra $A := \mathbb{C}\{a\} \simeq \mathbb{C}\{s\}$ of germs of holomorphic functions at the origin of $\mathbb{C}$. \\
We shall use the notation $ \mathbb{C}\{a\}$ when we consider $A$ as the sub-algebra of $End_c(\mathbb{C}\{s\})$ of $\mathbb{C}$-linear continuous endomorphisms of $\mathbb{C}\{s\}$ corresponding to multiplication by an element in $\mathbb{C}\{s\}$. So $a$ is the multiplication by $s$ acting on $\mathbb{C}\{s\}$.
\item The algebra $B := \mathbb{C}\{\{b\}\}$  is the closed  sub-algebra of the algebra $End_c(\mathbb{C}\{s\})$ of  $\mathbb{C}$-linear continuous  endomorphisms of $A$,  generated by $b$, the primitive without constant in $\mathbb{C}\{s\}$. So $b[f](s) := \int_0^s f(t)dt $.

\item The algebra $B[a]$ of polynomials in $a$ with coefficients in $B$, with the commutation relations $aS(b) - S(b)a = b^2S'(b)$ for each $S \in B$.\\
Note that we may write the coefficients in $B$ in the left or in the right without changing the algebra $B[a]$. The first choice corresponds to the left $B$-module structure on $B[a]$ and the second one to the left $\mathbb{C}[a]$-module structure of $B[a]$.\\

We shall use also the formal completions:
\item The algebra $\widehat{A} := \mathbb{C}[[a]]$,  the algebra\footnote{See Proposition 1 of Section 2 in \cite{[B.97]} for its relationship with $A_0[[b]]$.} $\widehat{B} := \mathbb{C}[[b]]$,  the algebra $\tilde{\mathcal{A}} := \widehat{B}[a]$ of polynomials in $a$ with coefficients in $\widehat{B}$, with the commutation relations $aS(b) - S(b)a = b^2S'(b)$ for each $S \in \widehat{B}$ and  the algebra $\widehat{\mathcal{A}}$ which the algebra of formal power series in $a$ and $b$ with the commutation relation $ab - ba = b^2$ (assuming that $a$ and $b$ are continuous for the topology defined by the valuation in (a, b)).
\end{enumerate}

 \parag{Two  easy facts} \begin{enumerate}
 \item The vector space of  homogenous polynomials in $(a,b)$ of degree $d$ in $\mathcal{A}_0$ admits either $a^{d-j}b^j, j \in [0,d]$ or $b^ja^{d-j}, j \in [0,d]$ as a basis.
\item  Any homogeneous polynomial in $(a, b)$  of degree $d$ which is  monic in $a$  in $\mathcal{A}_0$ may be written 
$$ P = (a - \lambda_1b)\dots (a - \lambda_db) $$
where $\lambda_1, \dots, \lambda_d$ are complex numbers (but such an expression in not unique, in general). See \cite{[B.09]}  Section 3.5 for the description of the various possibility to write a given $P$ in such a way.\\
 \end{enumerate}

Thanks to point $1.$ 	above we may define the numbers $\Gamma_{p,q}^j$ for $p, q \in \mathbb{N}$ and for $j \in [0, p]$  by the equation in $\mathcal{A}$
\begin{equation}
a^pb^q = \sum_{j=0}^p  \Gamma_{p,q}^jb^{q+j}a^{p-j} 
\end{equation}
It will be convenient to define $\Gamma_{p,q}^j = 0$ for $j > p$ or $j < 0$..

\begin{lemma}\label{11/9/21}
We have the following recursion relations:
\begin{align*}
& \Gamma_{p+1,q}^j = \Gamma_{p,q}^j + (q+j-1)\Gamma_{p,q}^{j-1} \quad \quad  \forall p, q \geq 0  \quad {\rm and} \quad  \forall j \in [1,p] \\
& {\rm and } \quad \Gamma_{p+1,q}^j = \Gamma_{p,q}^j + q\Gamma_{p,q+1}^{j-1} \quad \quad \forall p, q \geq 0  \quad  {\rm and} \quad  \forall j \in [1, p+1] .
\end{align*}
\end{lemma}

\parag{proof} To prove the first relation multiply  on the left the equality $(1)$ by $a$. This gives  
$\Gamma_{p+1, q}^j $ is the coefficient of $b^{q+j}a^{p+1-j}$ in the sum
\begin{align*}
& \sum_{h=0}^p \Gamma_{p,q}^hab^{q+h}a^{p-h} = \sum_{h=0}^p \Gamma_{p,q}^h(b^{q+h}a + (q+h)b^{q+h+1})a^{p-h} \\
& \qquad = \sum_{h=0}^p \Gamma_{p,q}^hb^{q+h}a^{p-h+1} +  \sum_{h=0}^p  \Gamma_{p,q}^h(q+h)b^{q+h+1})a^{p-h} 
\end{align*}
and we obtain the first relation.\\
For the second relation, multiply on the right the equality $(1)$ by $a$. This gives, since $a^pb^qa = a^{p+1}b^q + qa^pb^{q+1}$,
the result. $\hfill\blacksquare$\\
 
  Comparing the two relations  we obtain:
 $$ q\Gamma_{p,q+1}^{j-1} = (q+j-1)\Gamma_{p,q}^{j-1} \qquad \forall p \geq 1\quad \forall q \geq 0 \quad  \forall j \in [1, p+1] $$
 which leads to 
 $$ \Gamma_{p,q+1}^j = \frac{(q+j)!}{j!q!}\Gamma_{p,1}^j  \qquad \forall p \geq 1 \quad \forall q \geq 0 \quad  \forall j \in [0, p].$$
 As have
 $$ \Gamma_{p,0}^j = 1 \ {\rm for} \  j = 0 \quad {\rm and} \  \Gamma_{p,0}^j = 0 \  {\rm for} \  j \in [1,p] $$
 the first relation in the lemma above  gives
 $$ \Gamma_{p,1}^j = \frac{p!}{(p-j)!}$$
 So we have proved the following formula:
\begin{equation}
 \Gamma_{p,q}^j = \frac{(q+j-1)!}{j!(q-1)!}\frac{p!}{(p-j)!} \qquad \forall q \geq 1 \quad \forall p \geq 0 \quad {\rm and} \ \forall j \in [0, p]
 \end{equation}

 \parag{Remark} There is a $\mathbb{C}$-linear anti-automorphism $F$ of $\mathcal{A}$ defined by the following conditions
\begin{enumerate}
\item $F(xy) = F(y)F(x) \quad \forall  x, y \in \mathcal{A}$.
\item $F(a) = a, F(b) = -b \quad {\rm and} \quad F(1) = 1$
\end{enumerate}
Then apply $F$ to the relation $(1)$ gives
\begin{equation}
b^qa^p = \sum_{j=0}^p (-1)^j \Gamma_{p,q}^ja^{p-j}b^{q+j} 
\end{equation}
which inverts the relation $(1)$ when we consider it as the base change in the vector space of homogeneous polynomials of degree $m$  in (a,b) inside $\mathcal{A}$ for the basis $(a^pb^q)_{p+q = m}$ to the basis 
$(b^qa^p)_{p+q= m}$.

\begin{cor}\label{12/9/21}
For each $x \in \mathbb{C}$ we have the equality
\begin{equation}
(a + xb)^p = \sum_{j=0}^p \gamma_j(x)C_p^jb^ja^{p-j} 
\end{equation}
in the algebra $\mathcal{A}$, where $\gamma_j(x) := (x+j-1)(x+j-2) \cdots x $ for $j \in [0,p]$ (with $\gamma_0(x) \equiv 1$).
\end{cor}

\parag{proof} Since  both sides are  degree $p$ polynomials in $x$  with coefficients homogeneous of  degree $p$  in (a,b) in the $\mathbb{C}$-algebra $\mathcal{A}_0$, it is enough to prove  this formula for each $x = q \in \mathbb{N}$. Let us show first that we have the identity  $b^q(a+q b)^p = a^p b^q$ in $\mathcal{A}_0$:\\
This is an easy consequence of the fact that  $b(a+qb) = ab - b^2 + qb^2 = (a +(q-1)b)b$. \\
Then the formulas $(1), (2)$ and  $(3)$ imply $(4)$  for $x = q \in \mathbb{N}$,  because $b$ is not a zero divisor in $\mathcal{A}_0$. $\hfill \blacksquare$\\

Note that for $x = q \in \mathbb{N}^*$ we have $\gamma_j(q) = (q+j-1)!/(q-1)!$.

\parag{An example} As $\gamma_j(-1) = 0$ for $j \geq 2$ we obtain in $\mathcal{A}_0$ the relations:
$$ (a-b)^p = a^p - pba^{p-1} \quad \forall p \geq 2$$
This equality is easy to prove directly by induction on $p \in \mathbb{N}$, for instance using the identities $(a -b)^pb = ba^p$ and $a^pb = ba^p + pba^{p-1}b$.

\parag{Exercise} Show that for each integers $p, q \geq 0$ we have 
$$(a + xb)^pb^q = b^q(a + (x+q)b)^p.$$
 Then for $x = -q$ this gives $(a - qb)^pb^q = b^qa^p$.$\hfill \square$\\

\begin{prop}\label{produit}
Let $\Abc$ be the $\mathbb{C}-$vector space 
$$ \Abc := \{ \sum_{p,q} \gamma_{p,q}a^pb^q \/ \quad \exists R > 1, \ \exists C_R \quad {\rm such\ that}  \quad   \vert \gamma_{p,q}\vert \leq C_RR^{p+q}q! \}.$$
Assume that the variables $a$ and $b$ satisfy the commutation relation $ab - ba = b^2$. Then $\Abc$ is a sub-algebra of $\widehat{\mathcal{A}}$ of formal power series in the variables  $(a,b)$ with the relation $ab - ba = b^2$ (see above).\\
Moreover $\Abc$ is also described as the vector space:
$$ \Abc := \{ \sum_{p,q} \delta_{p,q}b^q a^p\/ \quad \exists R > 1, \  \exists D_R \quad {\rm such\ that} \quad   \vert \delta_{p,q}\vert \leq D_RR^{p+q}q! \}.$$
\end{prop}

\parag{Proof} Let $X := \sum_{p,q} \gamma_{p,q}a^pb^q $ and $Y := \sum_{p',q'} \delta_{p',q'}a^{p'}b^{q'} $ be in $\Abc$. Then the product $XY$ in $\widehat{\mathcal{A}}$ is given by
$$ XY = \sum_{m,n} \varepsilon_{m,n}a^mb^n $$
where
$$ \varepsilon_{m,n} = \sum_{p+p'-j = m, q+q'+j = n} (-1)^j \Gamma_{p',q}^j\gamma_{p,q}\delta_{p',q'} $$
using Formula $(3)$ to compute $b^qa^{p'}$.\\
 There exists positive constants   $R $ and $C_R$  large enough such that the following  estimates hold true
$$ \vert \gamma_{p,q}\vert \leq C_RR^{p+q}q! \quad {\rm and} \quad \vert \delta_{p,q}\vert \leq C'_RR^{p+q}q! $$
and then we obtain
\begin{align*}
&\vert  \varepsilon_{m,n} \vert \leq C_RC'_RR^{m+n}n! \sum_{p+p'-j = m, q+q'+j = n} \frac{(q+j-1)! p'! q! q'!}{n!(q-1)!j!(p'-j)!} \\
\end{align*}
In order to estimates the sum above, first note that  $q, q', j$  are at most equal to $n$ so that the triple $(q, q', j)$ takes at most $(n+1)^2$ values. Then $p$ and $p'$ are at most equal to $m+n$ and so the sum
has at most $(n+1)^2(m+n+1)$ terms. Each term is now bounded by
$$ q\frac{C_{p'}^j}{nC_{n-1}^{q'}} \leq C_{p'}^j \leq 2^{m+n} $$
using the fact that $q+j-1 = n-q'-1$. So there exists a constant $D_{3R}$ such that
$$ C_RC'_{R}(2R)^{m+n}(n+1)^2(m+n+1) \leq D_{3R}(3R)^{m+n} \quad \forall m,n \in \mathbb{N} $$
which implies 
 $$\vert  \varepsilon_{m,n} \vert  \leq D_{3R}(3R)^{m+n}n! $$
 showing that $XY$ is in $\Abc$.\\
 In order to write $X = \sum_{m,n} c_{m,n}b^na^m $ in $\widehat{\mathcal{A}}$ we use  Formula $(3) $ and we obtain
 $$ c_{m,n} = \sum_{p-j = m, q+j = n} (-1)^q \Gamma_{p,q}^j\gamma_{p,q} .$$
 The sum has at most $(n+1)^2 (m+n+1)$ terms and we have
 $$ \vert \Gamma_{p,q}^j\gamma_{p,q} \vert \leq  C_RR^{p+q}q! \frac{(q+j-1)!p!}{(q-1)!j!(p-j)!} \leq C_RR^{m+n} n! C_{p}^j \leq C_R(2R)^{m+n} n!$$
 and we conclude as above that there exists a constant $D_R$ such that  the following estimates hold true \quad
  $\vert c_{m,n}\vert \leq D_R (3R)^{m+n} n! \quad \forall (m,n) \in \mathbb{N}^2$.\\
  
  To prove the second assertion in the proposition take $X := \sum_{p, q} \gamma_{p,q} a^pb^q \in \Abc$ and write $X $ in  $\widehat{\mathcal{A}}$ as the series
  $$ \sum_{p,q} \delta_{p,q} b^qa^p .$$
  Then using Formula $(3)$ we have
  $$ \delta_{p,q} = \sum_{p',q',j}  \Gamma^j_{p', q'} \gamma_{p',q'} $$
  where the sum is on the triple  of non negative integers such that $(p', q', j)$ such that $ p = p' - j, q = q' + j$. So we have $p' + q' = p + q$ and  $j \in [0, q]$.
  Now using formula $(2)$ we have, assuming that $\vert \gamma_{p', q'}\vert \leq C_R R^{p'+q'} q'! $ :
  $$ \vert \delta_{p,q} \vert \leq C_R R^{p+q} q! \sum_{p', q', j}  \Gamma^j_{p', q'} \frac{q'!}{q!} .$$
  Now using Formula $(2)$ we obtain
  $$ \Gamma^j_{p', q'} \frac{q'!}{q!} =  \frac{q-j-1}{q} C_{p+j}^p \leq 2^{p+q}  $$
  and since our sum has at most $q+1$ terms, this gives 
  $$ \vert \delta_{p,q} \vert \leq C_R R^{p+q} (q+1)! $$
  completing the proof.  $\hfill \blacksquare$\\
  
  The following result is then immediate.
  
  \begin{cor}\label{20/3}
  The anti-automorphism $F$ of $\mathcal{A}$ (defined before Corollary \ref{12/9/21}) extends to a continuous  anti-automorphism $F_{conv}$ of $\Abc$ given by
  $$ F_{conv}(a^pb^q) = (-1)^qb^qa^p . \qquad  \qquad \qquad \qquad \qquad \qquad \qquad \blacksquare $$
  \end{cor} 
  
  This result allows to study only left modules on this algebra, since this anti-automorphism transform a right module in a left module.\\

  \parag{Remark} For $1 <  R < S$ given and any positive integers $m$ and $n$  there is a positive  constant $C_{m,n}$ such that
  $$ (\frac{R}{S})^{p+q}\frac{(q+m)!}{q!}(1+ p + q)^n \leq C_{m,n} \quad \forall (p,q) \in \mathbb{N}^2.$$
  So if $X := \sum_{p,q} x_{p,q}a^pb^q$ satisfies, for some given integers $m$ and $n$, the estimates 
  $$\vert x_{p,q}\vert \leq C_RR^{p+q}(q+m)! (1+p + q)^n \quad \forall (p,q) \in \mathbb{N}^2$$
   then $X$ is in $\Abc$.\\
  In particular, if  some $X \in \widehat{\mathcal{A}}$ \, is such that $Xb^m$ is in $\Abc$ then $X$ is also in $\Abc$. Moreover we may find $Y \in \Abc$ such that $Xb^m = b^m Y$.
   And the same is true if $b^mX$ is in $\Abc$, then $X$ is also in $\Abc$ and we may find $Z \in \Abc$ such that $b^mX = Zb^m$.\\
  The proof of the previous assertions is easy as it is enough to treat the case $m = 1$, using the second part of the previous theorem. Also,  if  some $X \in \widehat{\mathcal{A}}$ is such that  $Xa$ or $aX$ is  in $\Abc$ then $X$ is also in $\Abc$.$\hfill \square$\\
  
  \parag{Remark}  But contrary to $b\Abc = \Abc b$ the left and right ideals $a\Abc$ and $\Abc a$ are different: assume that $ab = ua$ for some $u \in \Abc$. Then we may replace $u$ by its homogeneous degree $1$  part in $(a, b)$. This implies that there exists $(\lambda, \mu) \in \mathbb{C}^2$ such that $ab =( \lambda a + \mu b) a $. Since $ba = (a-b)b$ this implies that
  $\lambda a^2$ is in $b\Abc = \Abc b$. This forces $\lambda = 0$ and then $ab = \mu ba $ which gives $(a - \mu a + b)b = 0$ and then $b = (\mu - 1)a$ which is impossible.$\hfill \square$\\

 \begin{lemma}\label{action 1}
 For any $X := \sum_{p,q} x_{p,q}a^p b^q \in \Abc$ define for $f(z) := \sum_{m=0}^\infty  t_mz^m$ in $\mathbb{C}\{z\}$
 $$X(f) := \sum_{m=0}^\infty u_mz^m \quad {\rm where} \quad  u_m := \sum_{p+q+r = m} \frac{r!}{(q+r)!} x_{p,q}t_r .$$
 Then $X(f)$ is in $\mathbb{C}\{z\}$.\\
 The corresponding  map of $\mathbb{C}$-algebras which sends $\Abc $ to the algebra of continuous endomorphisms of $\mathbb{C}\{z\}$,  is continuous and  injective.
 \end{lemma}
 
 \parag{proof} For  $R > 1$ large enough  there exists $C_R > 0$ and $D_R > 0$ such that 
 $$ \vert t_r\vert \leq C_RR^r \quad {\rm and} \quad \vert x_{p,q}\vert \leq D_RR^{p+q}q! .$$
 Then we obtain, since $a^pb^q[z^r] = \frac{r!q!}{(q+r)!} z^{p+q+r}$:
 $$ \vert u_m \vert \leq C_RD_R R^m \sum_{A_m}  \frac{r!q!}{(q+r)!} $$
 where the sum is taken on the set $A_m := \{(q, r) \in \mathbb{N}^2 \  / \  q+r \leq m \}$.
 The sum above has at most $(m+1)^2$ terms and each one is bounded by $1$ so for any $\rho > R$  there exists a constant $C_{\rho}$ such that $\vert u_m\vert \leq C_{\rho}\rho^m$, and this allows to conclude
 that $X(f)$ is in $\mathbb{C}\{z\}$.\\
 The continuity of this map follows from the estimates above.\\
 Consider now $X$ such that $X(z^r) = 0$ for each $r \in \mathbb{N}$.  The coefficient of $z^{m+r}$ in $X(z^r)$ is given by $P_m(z^r)$ where $P_m := \sum_{p+q= m} \gamma_{p,q}a^pb^q$ is the homogeneous part of degree $m$ in (a,b) of  $X$.  We know that  a non zero  homogeneous degree $m$ in (a,b) element in $\mathcal{A}$ may be written $b^j(a-\lambda_1b)\dots(a - \lambda_{m-j}b)$ where $\lambda_1, \dots, \lambda_{m-j}$ are complex numbers (see the "easy fact"  2 at the beginning of this section).\\
  Note that   $b$ is injective on $\mathbb{C}\{z\}$, and  that we have
 $$(a - \lambda b)(z^r) = (1 -\lambda/(r+1))z^{r+1}$$
 which vanishes only when $r+1 = \lambda$. So,  for any given $m$ there exists $r$ large enough such that $P_m(z^r) \not= 0$. Then, for $X \not= 0$ there exists  $m \in \mathbb{N}$ such that $P_m \not= 0$ and for $r$ large enough the coefficient of $z^{r+m}$ in $X(z^r)$ is not zero; so  the conclusion follows.$\hfill \blacksquare$\\
 
 \begin{lemma}\label{R bis} 
 For each non negative integer $N$ the vector space  $\mathcal{S}$ defined as the sub-set of series   $ X :=  \sum_{(p,q) \geq 0}  x_{p,q}a^pb^q $ in $\widehat{\mathcal{A}}$  satisfying
 $$ \exists R > 1, \  \exists C_R \quad {\rm such \ that} \quad \vert x_{p,q} \vert \leq C_R R^{p+q}\frac{(p+q)!}{p!}(1+p+q)^N \quad \forall p, q \geq 0 $$
 is equal to the algebra $\Abc$.
 \end{lemma}
 
 \parag{Proof} The inclusion  in $\Abc$ is clear because $p! q! \leq (p+q)! (1+p+q)^N$ for any integers $p, q, N \geq 0$. Conversely, if $X$ is in $\mathcal{S}$ we have
 $$\vert x_{p,q} \vert \leq C_R R^{p+q} C_{p+q}^q q! (1+p+q)^N \leq  \Gamma_NC_R (4R)^{p+q} q! $$
 because $C_{p+q}^q \leq 2^{p+q} $ and for each $N$ there exists $\Gamma_N$ large enough such that  $(1+p+q)^N \leq \Gamma_N 2^{p+q}$ for any $p$ and $q$ in 
 $\mathbb{N}$ .$\hfill \blacksquare$\\
 
 \begin{prop}\label{local}
 Let $X$ be in $ \Abc$ and assume that $x_{0,0} = 1$. Then $X$ is invertible in the algebra $\Abc$.
 \end{prop}
 
 \parag{Proof} Using the formula for the product in $\Abc$ obtained in the proof of Theorem \ref{produit} we see that, if $Y := \sum_{p',q'} y_{p',q'}a^{p'}b^{q'}$ is a formal inverse of $X$, then $Y$  is given by  $y_{0,0} = 1$ and for $(m,n) \not= (0,0)$ by
 $$y_{m,n} = \sum_{(j,p,q) \in A_{m,n}} (-1)^j \Gamma_{p',q}^j  x_{p,q}y_{p',q'}$$
 where  we define
  \begin{equation*}
   A_{m,n} = \{(p,q,j) \in \mathbb{N}^3 \ / \ \exists (p',q') \in B_{m,n} :  p+p'-j = m, \ q+q'+j = n \} 
   \end{equation*}
  where  $B_{m,n}$ is  the set of $(p',q') \in \mathbb{N}^2$ which are strictly smaller than $(m, n)$ for the lexicographical  order of $\mathbb{N}^2$.\\
  
   Choose $R > 1$ and $C_R > 0$ such that
  $$  \vert x_{p,q} \vert \leq C_R R^{p+q}\frac{(p+q)!}{p!} \quad\quad  \forall p, q \geq 0$$
    and assume that $S > R$ is large enough to satisfy satisfies
   $$C_R R/(S-R) \leq 1.$$
   We shall prove by induction on $(p',q') \in \mathbb{N}^2$, where we consider the lexicographical order on $\mathbb{N}^2$, the estimates
   \begin{equation*}
    \vert y_{p',q'} \vert \leq D_S S^{p'+q'} \frac{(p'+q')!}{p'!} \tag{$@$}
    \end{equation*}
   
     So assume that this estimates has been obtained for $(p',q') \in B_{m,n}$ for some constant $D_S$. Remark that for $m = n = 0$ we have only to ask that $D_S \geq 1$. \\
   Now fix $(m, n)$.  \\

  First consider the case where $n = 0$. So our induction assumption is that for any $m_1 < m$ and any $n_1 \geq 0$ we have already obtained the estimates $(@)$ for
  $(m_1, n_1)$ in $B_{m,0}$.
     Since we have
   $$ y_{m, 0} = \sum_{p=1}^m  \Gamma^0_{m-p, 0} x_{p, 0}y_{m-p, 0} $$ 
   we obtain, since $\Gamma^0_{p, 0} = 1$ and $RC_R/(S - R) \leq 1$ 
   
   $$ \vert y_{m,0} \vert \leq \sum_{p=1}^m C_RR^pD_SS^{m-p} \leq D_SS^m \sum_{p=1}^m C_R(R/S)^p \leq D_SS^m RC_R/(S -R) $$
   
   and then  $\vert y_{m,0}\vert \leq D_SS^m$.\\
   
   Now we have to consider the case $n \geq 1$ assuming that the estimates $(@)$ is already obtain for any  $(m_1, n_1)$ such that $m_1 < m$ and any $n_1\geq 0$ and for any  $(m, n_2)$ with $0 \leq n_2 < n$, so for $(m_1, n_1) \in B_{m,n}$.\\
    We   describe $A_{m,n}$  by fixing $j, q$ and $r := p+q$ so we have $p = r-q$,\\ $ p' = m+j-(r-q)$ and $q' = n-j-q$, and 
    $$  y_{m,n} = \sum_{(j,p,q) \in A_{m,n}} (-1)^j \Gamma_{p',q}^j  x_{p,q}y_{p',q'} $$
  $$ \sum_{p,q,j \in A_{m,n}} = \sum_{r=1}^{m+n} \sum_{q=0}^n\sum_{j = 0}^{m-r+q} .$$
   Then, since  $(p',q') \in B_{m,n}$ implies the estimates $(@)$,  we obtain
   
    \begin{align*}
  & \vert y_{m,n} \vert \leq D_SS^{m+n} \frac{(m+n)!}{m!} \sum_{p,q,j \in A_{m,n}}  C_R(R/S)^{p+q} \frac{(p+q)!(p'+q')!}{p!p'!} \frac{m!}{(m+n)!} \Gamma_{p',q}^j \\
  & \quad \leq D_SS^{m+n} \frac{(m+n)!}{m!} \big[ \sum_{r=1}^{m+n}  \frac{C_R(R/S)^r}{C_{m+n}^r}\Phi(r)\big] \tag{@@}
  \end{align*}
  where, fixing $r \in [1, m+n]$, we have to estimate, using that $\Gamma_{p',q}^j  = C_{q+j-1}^{q-1} \frac{p'!}{(p'-j)!} $, the quantity  
    \begin{equation*}
   \Phi(r) := \sum_{q=0}^n\sum_{j=0}^{n-q} \frac{m!}{p!p'!}\Gamma_{p',q}^j = \sum_{q=0}^n\sum_{j=0}^{n-q} \frac{m!}{p!p'!} C_{q+j-1}^{q-1} \frac{p'!}{(p'-j)!} .
  \end{equation*}
  Using the equality $m = p+p'-j$ we obtain
  \begin{equation*}
  \Phi(r)  = \sum_{q=0}^n\sum_{j=0}^{n-q} \frac{m!}{p!(p'-j)!}C_{q+j-1}^{q-1} =  \sum_{q=0}^n\sum_{j=0}^{n-q} C_m^{r-q}C_{q+j-1}^{q-1}. \tag{@@@}
  \end{equation*}
   
  Now we shall use the following elementary lemma.
  
  \begin{lemma}\label{combine}
  For any positive integers $ x, y $ we have
  $$ \sum_{j= 0}^{y} C_{x+j}^j = C_{x+y+1}^{x+1}.$$
  \end{lemma}
  
  \parag{Proof} As $C_{x+j}^j = C_{x+j}^x$ we are computing the coefficient of $a^x$ in the polynomial  $ \sum_{j=0}^{y} (1+a)^{x+j} $ which is equal to 
    $$(1+a)^{x}\frac{(1+a)^{y+1} - 1}{a} .$$
  So our sum is equal to the coefficient of $a^{x+1}$ in the polynomial $$ (1+a)^{x+y+1} - (1+a)^{x}$$
   and it is equal to $C_{x+y+1}^{x+1}$.$\hfill \blacksquare$\\
  
\parag{End of the induction of Proposition \ref{local}} So, thanks to the equality of Lemma \ref{combine}, for $x := q-1$ and  $y:= n-q$, which gives $\sum_{j=0}^{n-q} C_{q+j-1}^{q-1}  = C_n^q$, the sum in $(@@@)$ admit the estimates 
  $$ \Phi(r) \leq  \sum_{q=0}^n   C_m^{r-q}C_{n}^{q} \leq C_{m+n}^r  $$
 because $\sum_{p=0}^r C_m^pC_n^{r-p} = C_{m+n}^r $. Then the estimates $(@@)$ gives   
 $$\vert  y_{m,n} \vert \leq D_SS^{m+n} (m+n)! \frac{1}{m!} \sum_{r= 1}^{m+n}   \frac{C_R(R/S)^r}{C_{m+n}^r}\Phi(r) \leq D_SS^{m+n} (m+n)! \frac{1}{m!} $$
 as we choose $S$ in order that
 $$ C_R\sum_{r= 1}^{m+n} (R/S)^r \leq C_R\frac{R}{S - R} \leq 1$$
 we obtain that $Y$ is in $\Abc$. $\hfill \blacksquare$\\

 As a consequence, we see that $a\Abc+ b\Abc = \Abc b + \Abc a$ is the unique two sided maximal  ideal in $\Abc$ which is closed and with quotient
  $$\Abc\big/a\Abc + b\Abc \simeq \mathbb{C}.$$

 \begin{cor}\label{geometric 00} Let  $\Theta$ be a $(k,k)$-matrix with complex entries such that its spectrum is disjoint from $-\mathbb{N}$ and consider the simple pole  convergent (a,b)-module $\mathcal{E}(\Theta)$ introduced in  section 2.3 of \cite{[part I]}. So $ae = \Theta be$ where $e := (e_1, \dots, e_k)$ is a $B$-basis of $\mathcal{E}(\Theta)$.Then the continuous  actions of the algebras $\mathbb{C}[a]$ and $B$ on $\mathcal{E}(\Theta)$ are the restriction of a continuous action of $\Abc$ on $\mathcal{E}(\Theta)$.
 \end{cor}

 \parag{Proof} Since the right multiplication by an element of $B$ is a continuous endomorphism of left  $\Abc$-module on $\Abc$, it is enough to define the (left) action of $\Abc$ on the   $B$-basis 
$e_j, j \in [1,k]$ of $\mathcal{E}(\Theta)$:\\
 Indeed, if $X \in \Abc$ and $Z \in B^k$ the action of $X$ on $Ze$ is, by definition, the action of $XZ \in \Abc^k$ on $e$ where
 $XZ$ is defined by the left  action of $\Abc$ on $\Abc^k$  by left multiplication (using the fact that $B \subset \Abc$) and then we use the action of $\Abc^k$ on  each $e$  defined below.\\
 
So, using the second assertion of Proposition  \ref{produit},  let $X =  \sum_{p,q} x_{p,q} b^qa^p \in \Abc^k$; there exists constants  $R > 1$ and $C_R > 0$  with $\vert x_{p,q}\vert \leq C_R R^{p+q} q! \quad \forall (p,q) \in \mathbb{N}^2$.\\
 Since we have $a^pe = (\Theta +(p-1)Id)\circ \dots\circ \Theta b^pe$ we obtain 
 $$ Xe = \sum_{p,q} x_{p,q}b^qa^pe = \sum_{j=0}^\infty  z_jb^je $$
 where $z_j = \sum_{p+q = j} x_{p,q}(\Theta +(p-1)Id)\circ \dots\circ \Theta $. We have the estimate, for $\tau $ an integer larger than $\vert\vert \Theta \vert\vert$
 $$ \vert z_j \vert \leq C_RR^j j! \sum_{p=0}^j  \frac{(j-p)! (\tau + p-1)!}{j! (\tau -1)!} \leq C_RR^j j!\sum_{p=0}^j \frac{C_{\tau+p-1}^p}{C_j^p} $$
 and we obtain
 $$ \vert z_j \vert \leq  2^\tau C_R(2R)^j j! $$
  and then $Z := \sum_{j=0}^\infty  z_jb^j$ is in $B^{k}$ and $Xe : = Ze$ is well defined in the $B$-module with basis $e$.

  \parag{Remark}  The special case where, for $\alpha \in ]0,1] \cap \mathbb{Q}$ and $k \in \mathbb{N}^*$, we consider the $(k,k)$-matrix $\Theta_\alpha$ which gives $\mathcal{E}(\Theta_\alpha) = \Xi_\alpha^{(k-1)}$  is important since it shows, using the Embedding Theorem of \cite{[part I]} section 5, that geometric (a,b)-modules are canonically left $\Abc$-modules.\\

\subsection{The automorphisms $\tau_x$ of $\Abc$ and the Division Theorem}

\begin{prop}\label{26/3/25}
For any complex number $x \in \mathbb{C}$ there exists a unitary automorphism $\tau_x$ of the $\mathbb{C}$ algebra $\Abc$ which commutes with the left and right action
of $B$ on $\Abc$ defined by
$$ \tau_x(a) = (a + xb), \quad \tau_x(b) = b, \quad \tau_x(1) = 1 .$$
\end{prop}

\parag{Proof} So consider an element $X := \sum_{p, q}  \delta_{p, q} b^q a^p $ in $\Abc$. We want to prove that the element $Y \in \widehat{\mathcal{A}}$ defined by
$$ Y := \sum_{p, q}  \delta_{p, q} b^q(a + xb)^p $$
 is again in $\Abc$, where, thanks to Corollary \ref{12/9/21}, we have:
$$ (a + xb)^p =  \sum_{j=0}^p \gamma_j(x)C_{p}^jb^j a^{p-j}. $$
 So writing  $Y =  \sum_{p', q'}  \varepsilon_{p', q'} b^{q'}a^{p'} $ we have to estimates 
$$ \varepsilon_{p', q'} = \sum_{j = 0}^{q'} \delta_{p'+j, q'-j} \gamma_j(x)C_{p'+j}^j .$$
Since there exists $R > 1$ and $C_R > 0$ such that $\vert \delta_{p, q} \vert \leq C_R R^{p+q} q! $  we have  the estimates 
$ \vert \delta_{p'+j, q'-j} \vert \leq  C_R R^{p'+q'} (q' - j)!$ and then, assuming that $\vert x\vert \leq M + 1$ for some integer $M$ which gives the estimates 
 $$\vert \gamma_j(x)\vert \leq \frac{(M + j)!}{M!},$$
 we obtain:
$$ \vert \varepsilon_{p', q'} \vert  \leq C_R R^{p'+q'} q'! \sum_{j=0}^{q'}  \frac{(q'-j)!}{q'!}\frac{(M + j)!}{M!} \frac{(p'+j)!}{p'! j!} $$

Now we have
$$ \sum_{j=0}^{q'}  \frac{C_{M+j}^M C_{p'+j}^{p'}}{C_{q'}^j} \leq (q' +1)2^{M+q'}2^{p'+q'} \leq  2^M 5^{p'+q'}  \Gamma \qquad \forall (p', q') \in \mathbb{N}^2 $$
where $\Gamma$ is a universal  constant. This implies the estimates
$$  \vert \varepsilon_{p', q'} \vert  \leq  \Gamma C_{R} 2^M  (5R)^{p'+q'} q'! $$
and then $Y$ is in $\Abc$. \\
Since $(a+ xb)b - b(a + xb) = b^2$, the proof is complete. $\hfill \blacksquare$\\

First we shall work inside the algebra $\mathcal{A}_0$ of polynomials in the variables a and b with the commutation relation $ab - ba = b^2$.\\
Note that this algebra is integral and that any homogeneous element in $(a, b)$ of degree $m \geq 1$ monic "in a"  may be factorized as 
$$ P_m := (a-\lambda_1b)(a-\lambda_2b) \dots (a-\lambda_mb) $$ and  that the polynomial $\pi(x) := (x-(\lambda_1+ m -1))(x - (\lambda_2 + m +2))\dots (x - \lambda_m) $ satisfies the relation
$$ b^m\pi(b^{-1}a) = P_m $$
where the computation is made in the algebra $\mathcal{A}[b^{-1}]$.

\begin{lemma}\label{div.1}
For each  complex number $\lambda$ and each integer $m \in \mathbb{N}^*$ we have the equality
\begin{align*}
& \quad \quad a^m = Q_{m-1}(\lambda)(a - \lambda b) + R_m(\lambda) \quad  {\rm where} \quad \\
&  \quad \quad  Q_{m-1}(\lambda) = a^{m-1} + \lambda a^{m-2}b ± \dots + \lambda(\lambda+1) \dots (\lambda+ m-2)b^{m-1} \quad  {\rm and}  \\ 
&   \quad\quad R_m(\lambda) = \lambda(\lambda+1) \dots (\lambda+ m-1)b^m .
\end{align*}
\end{lemma}

\parag{Proof} For $m = 1$ the relation $a = (a-\lambda b) + \lambda b$ is clear and give $Q_0(\lambda) \equiv 1$ and $R_1(\lambda) = \lambda b$. So assume that the lemma is proved for $m \geq 1$ and multiply on the left by a. We obtain:
$$ a^{m+1} = aQ_{m-1}(\lambda) + aR_m(\lambda) = aQ_{m-1}(\lambda) + R_m(\lambda)(a+ m b). $$
Writing \quad  $R_m(\lambda)(a+m b) = R_m(\lambda)\big((a - \lambda b ) + (m + \lambda)b\big)$  \quad we obtain
$$ Q_m(\lambda) = aQ_{m-1}(\lambda) + R_m(\lambda) \quad {\rm and} \quad R_{m+1}(\lambda) = (\lambda + m) bR_m(\lambda) $$
completing the proof.$\hfill \blacksquare$

\begin{prop}\label{div.2}
Let $X := \sum_{(p,q) \in \mathbb{N}^2} x_{p,q}b^qa^p $ be an element in the algebra $\Abc$. Then for any  complex  number $\lambda$ there exists  a unique $Q$ in $\Abc$ and an unique  $R $ in $B$ such that the following equality holds in $\Abc$:
$$ X = Q(a - \lambda b) + R $$
\end{prop}

\parag{Proof}  To prove the uniqueness we have to show that $Q(a -\lambda b) + R = 0$ with $Q \in \Abc$ and $R \in  B$ implies $Q = 0$ and $R = 0$.\\
Assume that $Q$ is not $0$ and consider the minimal integer $m \geq 1$ such that  there is a non zero homogeneous term of degree $m-1$ in $(a,b)$  inside $Q$, denote it  $q_{m-1}$. Then the assumption implies  that the homogeneous term of degree $m$ in (a,b) inside $Q(a-\lambda b) + R$ vanishes. So we have
$ q_{m-1}(a - \lambda b) + r_mb^m = 0$ for some complex number $r_m$ which is the coefficient of $b^m$ in $R$.  For each integer $m$ the vector space of homogeneous elements of degree $d$ in $(a,b)$ admits the basis  
$$(a-\lambda b)^m, b(a -\lambda b)^{m-1}, \dots,  b^{m-1}(a -\lambda b), b^m$$ 
and the relation above gives a non trivial linear relation between these linearly independent  elements. This contradicts the assumption that $q_{m-1} \not= 0$. Then $Q = 0$ and so $R = 0$.\\
To prove the existence, it is enough to consider the case $\lambda = 0$ because if we may write
$$ \tau_\lambda(X) = Q a + R $$
with $Q \in \Abc$ and $R \in B$ then applying the automorphism $\tau_{-\lambda}$ gives
$$ X = \tau_{-\lambda}(Q)(a - \lambda b) + R $$
since $\tau_{-\lambda}$ acts on $B$ as the identity. \\
Now the case of the division by $a$ is obvious since we may write 
 $$X:=  \sum_{p \geq 0}  X_p(b)a^p = Q a + R $$
 with $Q :=  \sum_{p \geq 1}  X_p(b)a^{p-1} $ and $R := X_0(b)$ where $Q$ is clearly in $\Abc$.$\hfill \blacksquare$\\

\begin{thm}\label{div.3} Let $\lambda_1, \dots, \lambda_k$ be  complex  numbers and  let $S_1, \dots, S_{k}$ be invertible elements in $B$. Then define $P \in B[a]$ by
$$ P := (a - \lambda_1b)S_1(a -\lambda_2b)S_2 \dots (a - \lambda_kb)S_k .$$
Then for any $X$ in $\Abc$ there exists unique  $Q \in \Abc$ and $R \in B[a]$ with $deg_a(R) \leq k-1$ such that
$$ X = QP + R .$$
\end{thm}

\parag{Proof} The uniqueness is clear from the uniqueness statement in Proposition \ref{div.2} by an easy induction on $k$.\\
We shall prove the existence also  by induction on $k \geq 1$. For $k = 1$ it is enough to apply Proposition \ref{div.2} to $XS_1^{-1}$. So assume the proposition proved for $k-1$ and 
then write $$ X = Q_0(a - \lambda_2b)S_2 \dots (a-\lambda_kb)S_k + R_0$$
where $R_0$ has degree "in $a$" at most $k-2$. Now apply Proposition \ref{div.2}  for $(a - \lambda_1b)$ to $Q_0S_1^{-1}$. 
We obtain $$Q_0S_1^{-1} = Q_1(a - \lambda_1b) + R_1$$
where $R_1$ has degree $0$ in $a$. Then replacing $Q_0$ by $Q_1(a - \lambda_1b)S_1 + R_1S_1$ in the previous division we obtain
$$ X = Q_1(a- \lambda_1b)S_1(a -\lambda_2b)S_2 \dots (a - \lambda_kb)S_k + R_0 + R_1S_1(a -\lambda_2b)S_2\dots (a -\lambda_kb)S_k $$
this concludes the proof because
$$R := R_0 + R_1S_1(a -\lambda_2b)S_2\dots (a -\lambda_kb)S_k $$
has degree in $a$ at most equal to $k-1$.$\hfill \blacksquare$\\

\parag{Remark} As a direct consequence, we obtain that the quotient of $\Abc$ by the left ideal $\Abc P$ \ for such a $P$, coincides with the quotient $B[a]/B[a]P$ and is a free $B$-module with basis $1, a, \dots, a^{k-1}$. This is the general form of a  {\bf convergent frescos} with rank $k$ when the numbers $\lambda_j+j -k$ are rational and positive (see \cite{[part I]} and \cite{[part II]}) . \\
  So it will be equivalent, thanks to the previous results to consider geometric (a,b)-modules as left  $\Abc$-modules or as left $B[a]$-modules and  any $B[a]$-linear map between two geometric (a,b)-modules is $\Abc$-linear.
  
  \parag{An example} Let $\mathscr{A}$ be a finite subset in $]0, 1] \cap \mathbb{Q}$, $N$ be a non negative integer and $V$ a finite dimensional vector space.
  Let $\mathcal{E} := \Xi_\mathscr{A}^N \otimes_{\mathbb{C}} V$ the convergent (a,b)-module of  germs at $s =0$ convergent  multivalued series
  $$ \sum_{\alpha \in \mathscr{A}}\sum_{j=0}^{N} \sum_{m = 0}^\infty  s^{\alpha+m}(Log \, s)^j \otimes c_{\alpha, m, j} $$
  where the vector  $c_{\alpha, m, j}$ satisfy an estimates like $\vert\vert c_{\alpha, m, j}\vert\vert \leq C_R R^m $ for some $R > 1$ and some positive constant $C_R$, for some given norm on $V$. \\
  Then $\mathcal{E}$ is a free finite type $b$-module where the action of $B$ is defined as the "primitive without constant" (so  $B$ acts as $B := \mathbb{C}\{\{\partial_s^{-1}\}\}$
  see \cite{[part I]}). Thanks to the Division Theorem we may define the action of $\Abc$ on this convergent (a,b)-module. And since any geometric convergent (a,b)-module may be embedded as a $B[a]$-submodule of such an $\mathcal{E}$ (see \cite{[part I]} Theorem 5.1.4), this allows to consider any geometric convergent (a,b)-module as a  left $\Abc$-module.\\
  This fact will be used in the proof of Theorem \ref{5/10} (and also in Section 2.6) where such an $\mathcal{E}$ contains the integrals associated to $f : X \to D$ 
  $$ \int_{\gamma_s} \omega/df $$
  where $\omega \in \Omega^{p+1}$ satisfies $df\wedge  \omega = 0$ and $d\omega = 0$, where $(\gamma_s)_{s \in D¨*}$ is a horizontal family of  compact $p$-cycles in the fibers of $f$ and with $V := H^p(F, \mathbb{C})$ where $F$ is the Milnor fiber of $f$. The fact that such an integral is in such an $\mathcal{E}$ being consequence of the regularity of the Gauss-Manin connection and the positivity result of Malgrange (see \cite{[M.74]} , \cite{[M.75]} and also the Appendix of \cite{[B.84-b]})  .\\

  \subsection{Simple pole geometric differential systems and simple pole geometric (a,b)-modules}
  
  We generalize here to the geometric convergent case the Proposition 1.1 of \cite{[B.93]} which explains that any germ of a geometric  simple pole  linear differential system at the origin of $\mathbb{C}$ may be transformed in a simple pole geometric  convergent (a,b)-module.
  
  \begin{defn} \label{20/3}
   Fix an integer $p \geq 1$.  Let $M$ be a $(p, p)$ matrix with entries in $\mathbb{C}\{s\}$. We say that the germ at $s = 0$ of the  simple pole differential system
   \begin{equation}
  s \frac{dF}{ds} = M(s)F 
  \end{equation}
  is {\bf geometric} when the spectrum of the matrix $M(0)$ is in $\mathbb{Q}^{*+}$.
  \end{defn}
  
  \begin{thm}\label{3/25}
 Fix an integer $p \geq 1$.  Let $M$ be a $(p, p)$ matrix with entries in the algebra  $A := \mathbb{C}\{s\}$. Assume that the corresponding simple pole differential system is geometric.
  Let $\mathcal{E}$ be the free rank $p$ $B$-module with basis $e := (e_1, \dots, e_p)$. Then there exists an unique structure of simple pole geometric convergent (a,b)-module on  $\mathcal{E}$ such that we have
 $$ ae = M(a)be .$$
If this strucure is defined by the equality $ae = X(b)be$ where $X$ is a $(p,p)$ matrix with entries in $B := \mathbb{C}\{\{b\}\}$, then $X(0) = M(0)$ and the Bernstein polynomial of $\mathcal{E}$ is the minimal polynomial of $-M(0)$, the opposite of the residue matrix of the given  differential system.
 \end{thm}
  
 \parag{Proof}  Let $\varepsilon_1, \dots, \varepsilon_p$ be the horizontal basis of the differential system $(5)$. If $\mathscr{A}$ is the image in $\mathbb{Q} \cap ]0, 1]$ of the spectrum of $M(0)$, then each $\varepsilon_j$ is in $\Xi_\mathscr{A}^{(p)}$ and then consider the free rank $p$ $A$-sub-module of $\Xi_\mathscr{A}^{(p)}\otimes \mathbb{C}^p$ generated by $\frac{d\varepsilon_j}{ds}\otimes c_j$ where $c_1, \dots, c_p$ is the canonical basis of $\mathbb{C}^p$. Since $M(s)$ is invertible in $A \otimes End(\mathbb{C}^p)$, this $A$-sub-module is stable by $b $, the primitive without constant, and then is also a free rank $p$  $B$-sub-module. So it is a simple pole geometric (a,b)-module.$\hfill \blacksquare$\\
 
 \parag{Remark} It seems difficult to prove the complete generalization of the Proposition 1.1   of \cite{[B.93]}, so without the assumption on the spectrum of the matrix $M(0)$, because to estimate directly the coefficient of the formal solution $X \in \mathbb{C}[[b]]\times End(\mathbb{C}^p)$ given in {\it loc. cit.} seems quite hard\footnote{We know that without the spectral condition on $M(0)$ the formal solutions are convergent.}. In fact, this is not very surprising for the following reason:\\
 Assume that $e_1, \dots, e_p$ is the $B$-basis of a simple pole (a,b)-sub-module $\mathcal{E}$ of some  $\Xi_\mathscr{A}^{(N)}\otimes V$. Since  $a\mathcal{E} = b\mathcal{E}$, 
 $(be_1, \dots, be_p)$ is a $A$-basis of the free rank $p$ $A$-sub-module $a\mathcal{E}$. So there exists $M \in A \otimes End(\mathbb{C}^p)$ such that $ae = M(a)be$.\\
 But in practice, it is  complicated to relate the matrix $X \in B \otimes End(\mathbb{C}^p)$ such that $ae = X(b)be$ to the matrix $M$.\\
 
    \subsection{Some algebraic properties of $\Abc$}

We want to extend some algebraic properties of the algebra\footnote{The algebra $\tilde{\mathcal{A}}$  studied in \cite{[B.97]}.}  $(\mathbb{C}[a])[[b]]$ to $\Abc$. 

\parag{Terminology} We say that a left (resp right)  ideal $I \subset \Abc$  is {\bf normal} if $bx \in I$ (resp $xb \in I$) implies that $x \in I$.  

\begin{thm}\label{24/3}
The algebra $\Abc$ is left and right  noetherian. Moreover, any {\bf normal} left or right ideal is principal.
\end{thm}

\parag{Proof} So let consider first a normal left ideal $I$. Since $I \cap b\Abc = bI$, for any integer $q$ we have $b^q\Abc \cap I = b^qI$ and  the ideal $I$ is $b$-separated,
that is to say that
$$ \cap_{q \geq 0} (b^q\Abc \cap I) = \cap_{q \geq 0} b^q I = \{0\} $$
since $\cap_{q \geq 0} b^q \Abc = \{0\}$.\\
Consider the ideal  $J$ in $A := \mathbb{C}\{a\}$ of $T(a) \in A$ such that there exists $x \in I$ with $x = T(a)b^q + zb^{q+1}$  for some $q \in \mathbb{N}$. There exists
$p_0 \in \mathbb{N}$ such that $J  = (a^{p_0})$ and so there exists $x_0 \in I$ such that $x_0 = a^{p_0}b^{q_0} + zb^{q_0 +1} $ where $z $ is in $\Abc$. But since $I$ is assumed to be normal, we may choose\footnote{Write $x_0 = (a^{p_0} + zb)b^{q_0} = b^{q_0}(a^{p_0} + \zeta b)$ for some $\zeta \in \Abc$. Then $\xi_0 = a^{p_0} + \zeta b$ is in $I$ and we may replace $x_0$ by $\xi_0$.} $x_0$ such that  $q_0 = 0$.\\

Now consider any $y \in I \setminus bI$. After  left multiplication by an invertible element in $A$ we may write $y = a^{p_0 + p'_0} + bz$ with $z \in \Abc$ and $p'_0 \in \mathbb{N}$. Then $ y - a^{p'_0}x_0 =  by'_1$ and since $I$ is normal, we have $y_1 \in I$. Now write $y'_1 = b^{r_1}y_1$ with $q_0\in \mathbb{N}$ maximal.

\parag{Claim} There exist three sequences $(y_n)_{n \geq 0}$,  $(r_n)_{n \geq 1}$ and $(X_n)_{n\geq 0}$  with the following properties:
\begin{enumerate}
\item Each $y_n$ is in $I\setminus bI$ and $y_0 = y$ and with $y_1$ defined above.
\item Each $X_{n+1} - X_n$ is in $\Abc$  and has its valuation in  (a, b)  at least  equal to $n+1$.\\  Put  $X_0 = a^{p'_0} $.
\item Each $r_n$ is an integer such that $r_n \geq n$ and $r_1$  is defined above.
\item For each $n \geq 0$ we have  $y  = X_n x_0 + b^{r_{n+1}} y_{n+1}$.
\end{enumerate}

\parag{Proof of the claim} Remark first that $y_0 := y, y_1, r_1$ and $X_0$ are already constructed with the above properties. Assume the sequences $(X_n)_{n \geq 0}$, $(y_{n+1})_{n \geq 0}$ and $(r_{n+1})_{n \geq 0}$  are  constructed for integers at most $m$ for $m \geq 0$. \\
Then write $y_{m+1} = a^{p_0}T_m(a) + b z$ with $T_m \in A$. Then we have
$ y_{m+1} - T_m(a)x_0 = b\xi $ and so  $b^{r_{m+1}}y_{m+1}= b^{r_{m+1}}T_m(a)x_0 + b^{r_{m+1}+1}\xi $. This gives
$$ y = (X_m +  b^{r_{m+1}}T_m(a))x_0 +   b^{r_{m+1}+1}\theta .$$      
Then thanks to the normality of $I$  we have $\theta \in I$. Write $\theta = b^s y_{m+2}$ with $s \in \mathbb{N}$ and $y_{m+2} \in I \setminus bI$. 
Define $r_{m+2} = r_{m +1} + 1 +  s \geq m+2$ since $r_{m+1} \geq m+1$. We obtain    
$$ y = X_{m+1} x_0 + b^{r_{m+2}+1}y_{m+2} $$
where we  put $X_{m+1} := X_m + b^{r_{m+1}}T_m(a)$.  So our induction step, and then the claim is proved.$\hfill \blacksquare$\\

 So there exists  $U \in \widehat{\mathcal{A}}$:
 $$ U :=  X_0 +  \sum_{m=0}^\infty (X_{m+1} - X_m) = \lim_{m \to \infty}  X_m  $$
 
  such that  $Ux_0 = y$ in $ \widehat{\mathcal{A}}$,  since $ X_{m+1} - X_m$ is in $b^{m+1}\widehat{\mathcal{A}}$ for each $m$ 
  and then $I \subset \widehat{\mathcal{A}} x_0$. The next proposition allows to conclude that a normal left ideal is principal.

 \begin{prop}\label{28/3}
 Consider the element  $X := a^{p_0} + bz $ in $\Abc$  where $z $ is in $\Abc$. If $U \in \widehat{\mathcal{A}}$ is such that $Y := UX$ belongs to $\Abc$ then $U$ is also in 
 $\Abc$.
 \end{prop}
 
  \parag{proof} Put $X := \sum_{p, q} \gamma_{p, q} a^pb^q$ with $\gamma_{p, 0} = 0$ for $p \not= p_0$ and $\gamma_{p_0, 0} = 1$. There exists $R > 1$ and $C_R$ such that $\vert \gamma_{p, q}\vert \leq C_R R^{p+q} \frac{(p+q)! }{p!}$.\\
 Then write $U := \sum_{p, q} u_{p, q} a^pb^q $ and $UX = \sum_{m, n} \delta_{m, n} a^mb^n$ with the estimates (up to increase $R$ and $C_R$) 
 $$ \vert \delta_{m, n} \vert \leq C_R R^{m+n} \frac{(m+n)!}{m!}  .$$
 We shall show, by induction on $(m, n)$ using the lexicographic order on $\mathbb{N}^2$ that if there exist (for $R$, $C_R$ and  $p_0$ given)  positive constants $S$ and $D_S$ big enough to satisfying the conditions:
 \begin{itemize}
  \item  $S > 2R$ and $D_S \geq 2C_R(2R)^{p_0}$ and $D_S \geq  2RC_R\big/(S - R) $
  \end{itemize}
    and   such that for $(m, n) $ strictly less than $(m_0, n_0)$ for the lexicographical order  in $\mathbb{N}^2$ we have
 $$ \vert u_{m, n} \vert \leq D_S S^{m+n} \frac{(m+n)!}{m!} $$
 then  $$ \vert u_{m_0, n_0} \vert \leq D_S S^{m_0+n_0} \frac{(m_0+n_0)!}{n_0!} . $$
 
 The relation $Y = UX$ gives the equality
 
 \begin{equation*}
  \delta_{m, n} = \sum  u_{p, q}\gamma_{p', q'} a^pb^qa^{p'}b^{q'} = \sum_{j \leq p'}  u_{p, q}\gamma_{p', q'}\Gamma^j_{p', q}a^{p+p'-j}b^{q+q'+j}  \tag{@}
  \end{equation*}
 where the non negative integers  $p, q, p', q', j$ satisfy the relations
 \begin{equation*}  
 p+p' -j = m, \quad q + q' +j = n, \  0 \leq   j \leq p' \quad {\rm and} \  q' \geq  1 \quad{\rm when} \ p' \not= p_0. \tag{R}
 \end{equation*}
 So we have $q \leq n$ and also $p \leq m$ since $j \leq p'$.\\
 
 First consider the case where $n_0 = 0$. Then we want to estimate. $u_{m_0, 0}$. Consider the equality
 $$ \delta_{m_0+p_0, 0} = u_{m_0, 0} $$
 since $n_0 = 0$ implies that $q = q' = j = 0$ and we know that $\gamma_{p, 0} = 0$ for $p \not= p_0$. So the inequalities   $S > R, \quad D_S \geq C_R R^{p_0},$ which are obviously consequences of our requirements on $S$ and $D_S$,  are enough to conclude this case.\\
 Now assume that $n_0 \geq 1$ and that the estimates is proved for $(m, n)$ strictly less that $(m_0, n_0)$ (so either $m < m_0)$ or $m = m_0$ and $n < n_0$). The equality $(@)$ for $(m_0 + p_0, n_0)$  gives the estimates
 $$ \vert u_{m_0, n_0}\vert \leq \vert \delta_{m_0+p_0, n_0}\vert  + \vert \sum_{q' \geq 1, j \leq p'}   \Gamma^j_{p', q}u_{p, q}\gamma_{p', q'} \vert$$
 where in the sum $q' \geq 1 $ implies  $r' \geq 1$ (the case $q' = 0$ corresponds to $u_{m_0, n_0}$). So in each term of the sum above, $(p, q)$ is strictly less that $(m_0, n_0)$ for the lexicographical order and the inductive estimates applies. This gives
 $$ \vert u_{m_0, n_0}\vert \leq C_R R^{m_0+p_0+n_0} \frac{(m_0+p_0+n_0)!}{(m_0+p_0)!}+ \sum \Gamma^j_{p', q}D_S S^{p+q}\frac{(p+q)!}{p!} C_R R^{p'+q'} \frac{(p'+q')!}{p'!}.$$

  We   describe now the set  $A_{m_0,n_0}$  by fixing $j, q $ and $r := p+q$ as in Section 2.1. Denote by  $r' := m_0 + n_0 - r$ and since  $r' = p' + q'$ we have $r' \geq 1$. The  previous estimates give:
 \begin{align*}
 & \vert u_{m_0, n_0}\vert \leq D_S S^{m_0+n_0}\frac{(m_0+n_0)!}{m_0!} \Big[\frac{C_R R^{p_0}C_{m_0+p_0+n_0}^{m_0+p_0}}{D_S C_{m_0+n_0}^{m_0}}(\frac{R}{S})^{m_0+n_0} + 
  \sum_{r'=1 }^{m_0+n_0} \frac{C_R}{D_S}(\frac{R}{S})^{r'}\frac{ \Phi(r')}{C_{m_0+n_0}^r} \Big]\\
   \end{align*}
    where 
    
    $$\Phi(r') =  \sum_{q=0}^{n_0}\sum_{j=0}^{n_0-q} \frac{m_0!}{p!p'!} \Gamma_{p',q}^j = \sum_{q=0}^{n_0}\sum_{j=0}^{n_0-q} \frac{m_0!}{p!p'!} C_{q+j-1}^{q-1} \frac{p'!}{(p'-j)!}.$$
    
    We have, using Lemma \ref{combine}:
    
    $$ \Phi(r') =  \sum_{q=0}^{n_0}\sum_{j=0}^{n_0- q} C_{q+j-1}^j C_{m_0}^{r - q}  \leq  \sum_{q=0}^{n_0 -q}C_{n_0}^q C_{m_0}^{r - q} \leq C_{m_0+n_0}^r = C_{m_0+n_0}^{r'} $$
    
    Then  we obtain:
    $$\vert u_{m_0, n_0}\vert \leq D_S S^{m_0+n_0}\frac{(m_0+n_0)!}{m_0!}\Big[\frac{C_R R^{p_0}C_{m_0+p_0+n_0}^{m_0+p_0}}{D_S C_{m_0+n_0}^{m_0}}(\frac{R}{S})^{m_0+n_0}      +  \sum_{r'=1 }^{m_0+n_0} \frac{C_R}{D_S}(\frac{R}{S})^{r'}\Big.]$$
    Using now the fact that $\frac{C_{m_0+p_0+n_0}^{m_0+p_0}}{C_{m_0+n_0}^{m_0}} \leq 2^{m_0+n_0+p_0}$ we obtain:
    $$\vert u_{m_0, n_0}\vert \leq D_S S^{m_0+n_0}\frac{(m_0+n_0)!}{m_0!}\Big[ \frac{C_R}{D_S}(2R)^{p_0}\big(\frac{2R}{S}\big)^{m_0+n_0} +  \frac{C_RR}{D_S(S-R)}\Big].$$
    Since we choose $S > 2R$ and $D_S \geq 2C_R(2R)^{p_0}$ and $D_S \geq 2RC_R/(S - R)$, this concludes our induction.
    
              \parag{End of the proof of Theorem \ref{24/3}}   We already know that a normal left ideal is principal. Take now an arbitrary left ideal  $I$ and define
     \begin{align*}
     & J_p := \{ x \in \Abc \ / \ b^px \in I \} \quad {\rm for} \ p \in \mathbb{N} \quad {\rm and} \\
     & J_\infty := \cup_{p \geq 0} \ J_p .
     \end{align*}
     Then each $J_p$ is a left ideal in $\Abc$ because if $b^px \in I$ then $b^pyx = zb^px \in I$ since for any $y \in \Abc$ and any $p \in \mathbb{N}$ there exists $z \in \Abc$ such that $zb^p = b^py$. So we obtain an increasing sequence of left ideals $I \subset J_1 \subset \cdots \subset J_p \subset J_{p+1} \subset \cdots \subset J_\infty $.\\
     Moreover $J_\infty$ is normal, since $bx \in J_\infty$ implies $bx \in J_p$ for some $p$ and then $b^{p+1}x \in I$, so $x \in J_{p+1} \subset J_\infty$.\\
     Let $u$ be a generator of $J_\infty$. Then there exists $q \in \mathbb{N}$ such that $u $ is in $ J_q$ and so $J_q = J_\infty$.\\
     Now we have $b^qu \in I$ and so $b^qJ_\infty \subset I \subset J_\infty \simeq \Abc u$. So the quotient $I\big/b^q\Abc u$ is a sub-module of $J_\infty\big/b^qJ_\infty$. But this last left $\Abc$-module is isomorphic to $\Abc \big/ b^q\Abc$  as a $A$-module. And this $A$-module  is isomorphic to $\oplus_{h=0}^{q-1} Ab^h$ which is a free finite type $A$-module. Then $I\big/b^q\Abc u$ is a finite type $A$-module, since $A = \mathbb{C}\{a\}$ is noetherian.\\
      Now if $g_1, \dots, g_m$ is a generator system of this $A$-module then $g_1, \dots, g_m, b^qu$ is a generator of the ideal $I$:\\
     Let $z \in I$. Then the image of $z$ in $I\big/b^q\Abc u$ may be written $[z] = \sum_{j=1}^m t_jg_j $  for some $t_j \in A$ and this means that $z - \sum_{j=1}^m t_jg_j$ is in $b^q\Abc u = \Abc b^qu$ and then
     $z = \sum_{j=1}^m t_jg_j + \theta b^qu$, proving our claim.
     So $\Abc$ is left noetherian. \\
      The anti-automorphism $F$ (see Corollary \ref{20/3}) gives the right noetherianity of $\Abc$.$\hfill \blacksquare$\\

\section{Brieskorn Modules and Vanishing Cycles}

 \subsection{Preliminaries.}
 
First  we generalize to the convergent case the results of  Section 2 of \cite{[B.06]}. \\   
 The situation we  consider is the following : let \ $M$ \ be a connected  complex manifold of dimension \ $n +1, n \geq 1$ \ and \ $f : M \to \mathbb{C}$ \ a non constant holomorphic function.  We denote by $S$ the subset $\{df = 0\}$ and by  $Y$ the subset  $ f^{-1}(0)$ and we assume that $S \subset Y$. Moreover we assume that $f = 0$ is a reduced equation of $Y$\\
 
We introduce the following complexes of sheaves supported by \ $Y $. 

 \begin{enumerate}
 \item  The topological restriction \ $(\tilde{\Omega}^{\bullet}, d^{\bullet})$ \ to $Y$ of the usual holomorphic de Rham complex on  \ $M$.
 \item The sub-complexes \ $(K^{\bullet}, d^{\bullet})$ \ and \ $(I^{\bullet}, d^{\bullet})$ \ of \  $(\tilde{\Omega}^{\bullet}, d^{\bullet})$ \  where the sub-sheaves \ $K^p$ \ and \ $I^{p+1}$ \ are defined for each \ $p \in \mathbb{N}$ \  respectively as the kernel and the image of the map
 $$  \wedge df : \tilde{\Omega}^p \to \tilde{\Omega}^{p+1} $$
 given par exterior multiplication by \ $df$.  We have the exact sequence of complexes:
 \begin{equation*} 0 \to (K^{\bullet}, d^{\bullet}) \to (\tilde{\Omega}^{\bullet},  d^{\bullet}) \to (I^{\bullet}, d^{\bullet})[+1]  \to 0. \tag{1}
 \end{equation*}
 Note that \ $K^0$ \ and \ $I^0$ \ are zero by definition.
 \item The natural inclusions \ $I^p \subset K^p$ \ for all \ $p \geq 0$ \ are compatible with the differential \ $d$. This leads to an exact sequence of complexes
 \begin{equation*}
 0 \to (I^{\bullet}, d^{\bullet}) \to (K^{\bullet}, d^{\bullet}) \to ([K/I]^{\bullet}, d^{\bullet}) \to 0 .\tag{2}
 \end{equation*}
 \item We have a natural inclusion \ $f^{-1}(\tilde{\Omega}_{\mathbb{C}}^1) \subset K^1\cap Ker\, d$,    where $\tilde{\Omega}_{\mathbb{C}}^1$ is the topological restriction to $\{0\}$ of the sheaf $\Omega_{\mathbb{C}}^1$, and this gives a sub-complex (concentrated in degree $1$ ) of \ $(K^{\bullet}, d^{\bullet})$. As in \cite{[B.06]}, we shall consider also the  quotient complex \ $(\tilde{K}^{\bullet}, d^{\bullet})$ \  defined by  the exact sequence of complexes:
 \begin{equation*}
  0 \to f^{-1}(\tilde{\Omega}_{\mathbb{C}}^1) \to (K^{\bullet}, d^{\bullet}) \to (\tilde{K}^{\bullet}, d^{\bullet}) \to 0 . \tag{3}
  \end{equation*}

 We make the assumption here that \ $f = 0 $ \ is a reduced equation of \ $Y$,  so the cohomology sheaf in degree 1  of the complex \ $(K^{\bullet}, d^{\bullet})$, which is equal to \ $K^1 \cap Ker\, d$,  co{\'i}ncides with \ $f^{-1}(\tilde{\Omega}_{\mathbb{C}}^1)$. So $\tilde{K}^1 = 0$ and  the complex \ $ (\tilde{K}^{\bullet}, d^{\bullet})$ \  has a  zero cohomology sheaf in degree 1.
  \end{enumerate} 
  Recall now that we have on the cohomology sheaves of the following complexes \\
   $(K^{\bullet}, d^{\bullet}), (I^{\bullet}, d^{\bullet}), ([K/I]^{\bullet}, d^{\bullet}) $ \ and \ $f^{-1}(\tilde{\Omega}_{\mathbb{C}}^1), (\tilde{K}^{\bullet}, d^{\bullet})$ \ natural operations \ $a$ \ and \ $b$ \ with the relation \ $ab - ba = b^2$. They are defined in a na{\"i}ve way by 
  $$  a : = \times f \quad {\rm and} \quad  b : = \wedge df \circ d^{-1} .$$
  The definition of \ $''a''$ \ makes sens obviously. Let me precise the definition of \ $''b''$ \ first in  the case of the $p$-th cohomology sheaf  \ $\mathcal{H}^p(K^{\bullet}, d^{\bullet})$ \ with \ $p \geq 2$:\\ if \ $x \in K^p \cap Ker\, d$ \ write \ $x = d\xi$ \ with \ $\xi \in \tilde{\Omega}^{p-1}$ \ and let \ $b[x] : = [df\wedge \xi]$. The reader will check easily that this makes sens.\\
  For \ $p = 1$ \ we shall choose \ $\xi \in \tilde{\Omega}^0$ \ such that \ $\xi = 0$ \ on the smooth part of  \ $Y$. This is possible because the condition \ $df \wedge d\xi = 0 $ \ allows such a choice : near a smooth point of \ $Y$ \ we can choose coordinates such \ $ f = x_0$ \ and the condition on \ $\xi$ \ means independence of \ $x_1, \cdots, x_n$. Then \ $\xi$ \ has to be  locally constant on \ $Y$ \ which is locally connected. So we may kill the value of such a germ \ $\xi$ \ along \ $Y$.\\
  The case of the cohomology sheaf of  the complex \ $(I^{\bullet}, d^{\bullet})$ \  is reduced to the previous one using the next  lemma.
  
  \begin{lemma}\label{tilde b}
  For each \ $p \geq 0$ \ there is a natural injective map
  $$\tilde{b} :  \mathcal{H}^p(K^{\bullet}, d^{\bullet}) \to \mathcal{H}^p(I^{\bullet}, d^{\bullet})$$
  which satisfies the relation \ $a \tilde{b} = \tilde{b}(b + a) $. For \ $p \geq  2$ \ this map is bijective.
  \end{lemma}
  
  \parag{Proof} Let \ $x \in K^p \cap Ker\, d $ \ and write \ $x = d\xi $ \ where \ $x \in \Omega^{p-1}$ \ (with \ $\xi = 0$ \ on \ $Y$ \ if \ $p = 1$), and set \ $\tilde{b}([x]) : = [df\wedge \xi] \in \mathcal{H}^p(I^{\bullet}, d^{\bullet})$. This is independent on the choice of \ $\xi$ \ because, for \ $p \geq 2$, adding \ $d\eta$ \ to \ $\xi$ \ does not modify the result as \ $[df\wedge d\eta] = 0 $. For \ $p =1$ \ remark that our choice of \ $\xi$ \ is unique.\\
   This is also independent of the the choice of \ $x $ \ in \ $ [x] \in \mathcal{H}^p(K^{\bullet}, d^{\bullet})$ \  because adding \ $\theta \in K^{p-1}$ \ to \ $\xi$ \ does not change \ $df \wedge \xi$.\\
   Assume \ $\tilde{b}([x]) = 0 $ \ in \ $ \mathcal{H}^p(I^{\bullet}, d^{\bullet})$ ;  this means that for $p \geq 2$  we may find \ $\alpha \in \Omega^{p-2}$ \ such \ $df \wedge \xi = df \wedge d\alpha$ (and   $df \wedge \xi = 0$ and so $\xi = 0$ for $p = 1$). But then, \ $\xi - d\alpha $ \ lies in \ $K^{p-1}$ \ and \ $x = d(\xi - d\alpha ) $ \ shows that \ $[x] = 0$. So \ $\tilde{b}$ \ is injective.\\
  Assume now \ $p \geq 2$.  If \ $df\wedge \eta $ \ is in \ $I^p \cap Ker\, d$, then \ $df \wedge d\eta = 0 $ \ and \ $y : = d\eta $ \ lies in \ $K^p \cap Ker\, d$ \ and defines a class \ $[y] \in  \mathcal{H}^p(K^{\bullet}, d^{\bullet}) $ \ whose image by \ $\tilde{b}$ \ is \ $[df\wedge \eta] $. This shows the surjectivity of \ $\tilde{b}$ \ for \ $p \geq 2$.\\
   For \ $p=1$ \ the map \ $\tilde{b}$ \ is not surjective (see the remark below). The co-image of $\tilde{b}$ for $p = 1$ is isomorphic to $\mathbb{C}_Y df $.\\
  To finish the proof let us  to compute \ $\tilde{b}(a[x] + b[x])$. Writing again \ $x = d\xi$, we get
   $$ a[x] + b[x] =[ fd\xi + df \wedge \xi] = [d(f\xi)] $$
   and so
   $$ \tilde{b}( a[x] + b[x] ) = [ df \wedge f\xi ] = a\tilde{b}([x]) $$
   which concludes the proof. $\hfill \blacksquare$
   
   \bigskip
   
   Denote by \ $i :  (I^{\bullet}, d^{\bullet}) \to (K^{\bullet}, d^{\bullet})$ \ the natural inclusion and define the action of \ $b$ \ on \ $\mathcal{H}^p(I^{\bullet}, d^{\bullet})$ \ by \ $b : = \tilde{b}\circ \mathcal{H}^p(i) $. As \ $i$ \ is \ $a-$linear, we deduce the relation \ $ab - ba = b^2$ \ on \ $\mathcal{H}^p(I^{\bullet}, d^{\bullet})$ \ from the relation of the previous lemma. \\
   
   The action of \ $a$ \ on the complex \ $ ([K/I]^{\bullet}, d^{\bullet}) $ \ is obvious and the action of \ $b$ \ is zero.\\
   
   The action of \ $a$ \ and \ $b$ \ on \ $f^{-1}(\tilde{\Omega}_{\mathbb{C}}^1) \simeq \mathcal{E}_1\otimes \mathbb{C}_Y$ \ are the obvious one, where \ $\mathcal{E}_1$ \ is the rank 1 (convergent)  (a,b)-module with generator \ $e_1$ \ satisfying \ $a.e_1 = b.e_1$ \ (or, if you prefer, \ $\mathcal{E}_1 : = \Abc\big/\Abc(a - b)$). \\
   Remark that the natural inclusion \ $f^{-1}(\tilde{\Omega}^1_{\mathbb{C}}) \hookrightarrow (K^{\bullet}, d^{\bullet})$ \ is compatible with the actions of \ $a$ \ and \ $b$. The actions of \ $a$ \ and \ $b$ \ on \ $\mathcal{H}^1(\tilde{K}^{\bullet}, d^{\bullet}) $ \ are simply  induced by the corresponding actions on \ $\mathcal{H}^1(K^{\bullet}, d^{\bullet})$.
   
   \parag{Remark} The exact sequence of complexes (2) induces  for any \ $p \geq 2$ \  a bijection
   $$ \partial^p : \mathcal{H}^p(I^{\bullet}, d^{\bullet}) \to \mathcal{H}^p(K^{\bullet}, d^{\bullet})$$
   and a short exact sequence 
   \begin{equation*}
    0 \to \mathbb{C}_Y df \to \mathcal{H}^1(I^{\bullet}, d^{\bullet}) \overset{\partial^1}{\to} \mathcal{H}^1(K^{\bullet}, d^{\bullet}) \to 0 \tag{@}
    \end{equation*}
   because of the (holomorphic) De Rham Lemma. Let us check  that for \ $p \geq 2$ \ we have \ $\partial^p = (\tilde{b})^{-1}$ \ and that for \ $p =1$ \ we have \ $\partial^1\circ \tilde{b} = Id$. If \ $x = d\xi \in K^p \cap Ker\, d$ \ then \ $\tilde{b}([x]) = [df\wedge \xi]$ \ and \ $\partial^p[df\wedge\xi] = [d\xi]$. So \ $ \partial^p\circ\tilde{b} = Id \quad \forall p \geq 0$. For \ $p \geq 2$ \ and \ $df\wedge\alpha \in I^p \cap Ker\, d$ \ we have \ $\partial^p[df\wedge\alpha] = [d\alpha]$ \ and \ $\tilde{b}[d\alpha] = [df\wedge\alpha]$, so \ $\tilde{b}\circ\partial^p = Id$. For \ $p = 1$ \ we have \ $\tilde{b}[d\alpha] = [df\wedge(\alpha - \alpha_0)]$ \ where \ $\alpha_0 \in \mathbb{C}$ \ is such that \ $\alpha_Y = \alpha_0$.  This shows that in degree 1 \ $\tilde{b}$ \ gives a canonical splitting of the exact sequence \ $(@)$.

\subsection{Action of $\Abc$ on $\mathcal{O}_{M\vert Y}\{\{b\}\}$}

Let $M$ be a complex manifold and $f : M \to \mathbb{C}$ be a holomorphic function. Denote by $Y := \{ f = 0 \}$ the zero set of $f$.  As our construction is local on $Y$, we may assume that $M $ is an open set in $\mathbb{C}^{n+1}$ and replace $M$ and $f$  by a Milnor representative $f : X \to D$ of $f$ at a given point $y_0 = 0$ in $\mathbb{C}^{n+1}$.

\begin{defn}\label{28/9}
We define the sheaf $\mathcal{O}_{M\vert Y}\{\{b\}\}$ as follows:\\
A germ of section of this sheaf at the point $y_0 \in Y$ is the set of all  series like  $g := \sum_{m = 0}^\infty b^mg_m(x)$ where the $g_m, m \geq 0,$ are holomorphic functions such that there exists   a compact neighborhood $K$  of $y_0$ in $M$ on which each $g_m$ is  continuous and holomorphic in its interior points and such that there exists $R > 1$ and $C_R > 0$  for which  the following estimates estimates holds true
\begin{equation*}
\vert\vert g_m \vert\vert_K \leq C_RR^m m!  \tag{E1}
\end{equation*} 
\end{defn}
Here $\vert\vert g_m \vert\vert_K$ denotes the sup-norm on $K$ of the algebra $H(K)$ of the continuous functions on $K$ which are holomorphic in the interior points of $K$.\\
Note that this defines a sheaf on $Y$ since the existence of $K,R$ and $ C_R$ implies that such a series defines also a germ of section of $\mathcal{O}_{M\vert Y}\{\{b\}\}$ at each point $y $ in the interior of $K\cap Y$.\\
Compare with the definition of the sheaf $\Omega^\bullet_{M\vert Y}[[b]]$ given in \cite{[B.08]} whose germ at $y_0$ is 
$$\lim_{\underset{y_0 \in U}{\longrightarrow}} \Omega_M(U)[[b]]$$
which does not coincides with $\Omega_{M, y_0}[[b]]$ in general.

\begin{lemma}\label{easy 1}
The sheaf $\mathcal{O}_{M\vert Y}\{\{b\}\}$ is a sheaf of left $B:= \mathbb{C}\{\{b\}\}$-modules by the obvious rule $bg := \sum_{m=0}^\infty b^{m+1}g_m(x)$.
\end{lemma}

\parag{Proof} Let $T(b) = \sum_{q=0}^\infty c_qb^q \in B$. So there exists $S > 1$ and $\Gamma_S > 0$ such that $\vert c_q\vert \leq \Gamma_SS^q q! $ \  for each $q \geq 0$.
Assume that $g$ is as in the definition above and assume that $S \geq R$ (we may  increase the choice of $S$). Then we have
$$ T[g] = \sum_{n=0}^\infty b^n \Big(\sum_{q+m = n} c_qg_m(x)\Big) $$
and the following estimates hold true:
\begin{align*}
& \big\vert\big\vert \sum_{q = 0}^n c_q g_{n-q}\big\vert\big\vert_K \leq \sum_{q=0}^n \Gamma_S S^q q! C_R R^{n-q} (n-q)! \\
& \big\vert\big\vert \sum_{q = 0}^n c_q g_{n-q}\big\vert\big\vert_K \leq \Gamma_SC_RS^n n!\Big(\sum_{q=0}^n \frac{q!(n-q)!}{n!}\big(\frac{R}{S}\big)^{n-q}\Big) \\
& \big\vert\big\vert \sum_{q = 0}^n c_q g_{n-q}\big\vert\big\vert_K \leq  C_R\Gamma_S\frac{S}{S-R} S^n n!
\end{align*}
and this shows that $T[g]$ is again a section on the interior of $K$ of the sheaf  $\mathcal{O}_{M\vert Y}\{\{b\}\}$.$\hfill\blacksquare$\\

We have a natural action of $a$ on the sheaf  $\mathcal{O}_{M\vert Y}\{\{b\}\}$ defined by the product by $f$, taking in account that we want that $ab - ba = b^2$. So it is defined  by $$ab^mg_m(x) := b^mf(x)g_m(x) + mb^{m+1}g_m(x),$$
 and  using the bound $\vert\vert f\vert\vert_K := F$  we obtain, using the equality
$$ ag =  f(x)g_0(x) + \sum_{m=1}^\infty b^m\big(f(x)g_m(x) + (m-1)g_{m-1}(x)\big) $$
and  the existence of $R > 1$ and of $C_R > 0$ such that $\vert\vert g_m\vert\vert_K \leq C_R R^m m!$, 
the estimates $\vert\vert (ag)_0\vert\vert_K \leq FC_R$ and , for  each $m \geq 1$, the estimates:
$$\vert\vert  (ag)_m\vert\vert_K = \vert\vert f_m + (m-1)g_{m-1}\vert\vert_K \leq (F + 1) C_RR^m m!  .$$

The next proposition shows that this leads to a left  action of the algebra $\Abc$ on the sheaf  $\mathcal{O}_{M\vert Y}\{\{b\}\}$  which is compatible with the left action of $B$ given by Lemma \ref{28/9}.

\begin{prop}\label{Action}
The sheaf $\mathcal{O}_{M\vert Y}\{\{b\}\}$ is a left $\Abc$-module where $B$ acts by the left  action defined in Lemma  \ref{28/9}. and where the action of $a$ is defined by the formula 
\begin{equation*}
ab^mg_m(x) := b^mf(x)g_m(x) + mb^{m+1}g_m(x) \tag{E2}
\end{equation*}
\end{prop}

\parag{Proof} So consider $g = \sum_{m=0}^\infty b^mg_m(x)$ a section of the sheaf $\mathcal{O}_M\{\{b\}\}_{\vert Y}$ and $T := \sum_{p, q} c_{p, q}a^pb^q$ an element in $\Abc$. So there exist  $R > 1$ and $C_R > 0$ such that 
$$ \vert c_{p, q}\vert \leq C_R R^{p+q} q! $$
We may find  $S$ large enough and $\Gamma_S > 0$ such that

$$  S > 2R \quad {\rm and} \quad \vert\vert g_m \vert\vert_K \leq \Gamma_S S^m m! $$
Since we may shrink $K$ around $y_0 \in Y$ we may assume moreover  that 
$$\vert\vert f\vert\vert_K \leq S^{-1}.$$
To compute $a^pb^qb^mg_m(x)$ we use the formula 
$$ a^pb^{q+m} = \sum_{j = 0}^p \Gamma^j_{p, q+m} b^{q+m+j}a^{p-j}$$
and the fact that $a^{p-j}g_m(x) = f^{p-j}(x)g_m(x) $. Put $T[g] := \sum_{n=0}^\infty  b^nt_n(x) $. Then we have, with $r := p-j$
$$ t_n(x) = \sum_{m+q+j = n} \sum_{r=0}^\infty  c_{r+j, q}\Gamma^j_{r+j, q+m}f^{r}(x)g_m(x) $$
	and so, with $m = n-q-j$  we obtain:
\begin{align*}
& \vert\vert t_n\vert\vert_K \leq \sum_{m+q+j = n} \sum_{r= 0}^\infty  C_RR^{r+q+j}q! C_{r+j}^j\frac{(n-1)!}{(n-j-1)!}S^{-r}\Gamma_S S^{n-q-j}(n-q-j)! \\
&  \vert\vert t_n\vert\vert_K \leq C_R\Gamma_S S^n n! \Big[\sum_{q+j = 0}^n  \big(\frac{n-j}{n}\big)\frac{q!(n-q-j)!}{(n-j)!}\big(\frac{R}{S}\big)^{q+j} \Big(\sum_{r=0}^\infty  C_{r+j}^j\big(\frac{R}{S}\big)^r \Big)\Big]
\end{align*}
Now, using\footnote{ Exercise: Prove that for $x \in [0,1[$  the equality $\sum_{r=0}^\infty C^j_{r+j} x^r = (\frac{1}{1-x})^{j+1}$ holds true for any $j \in \mathbb{N}$.} 

$$\big(\frac{n-j}{n}\big) \leq 1, \quad C_{n-j}^q \geq 1, \quad {\rm and} \quad  \sum_{r=0}^\infty C^j_{r+j}\big(\frac{R}{S}\big)^r = \big(\frac{S}{S - R}\big)^{j+1} $$ 
we get 
$$ \frac{S}{S - R}\sum_{q+j = n} \big(\frac{R}{S - R}\big)^j\big(\frac{R}{S}\big)^q  \leq  \frac{S}{S - R}\sum_{s = 0}^\infty \big(\frac{R}{S-R}\big)^s = \frac{S}{S-2R}. $$
So we obtain
$$  \vert\vert t_n\vert\vert_K \leq C_R\Gamma_S\frac{R}{S-2R} S^n n! $$
This implies that $T[g]$ is a section on the interior of $K$ of the sheaf $\mathcal{O}_{M\vert Y}\{\{b\}\}$.$\hfill \blacksquare$\\

We now define  on $Y$ the sheaf of $\Abc$-modules $\Omega_{M\vert Y}^p\{\{b\}\}$ for each integer $p \geq 0$. The action of $\Abc$ on
$$\Omega_{M\vert Y}^p\{\{b\}\} :=  \mathcal{O}_{M\vert Y}\{\{b\}\} \otimes_{\mathcal{O}_{M\vert Y}} \Omega_{M\vert Y}^p $$
 is defined  by action of $\Abc$ on the left part of  the tensor product

So a germ of  section $\omega$ at the point $y_0 \in Y$  of this sheaf is a series  of the type $\omega := \sum_{m=0}^\infty b^m\otimes \omega_m$ where each $\omega_m$ is in the Banach space $H(K, \Omega^p_M)$ for some given  compact neighborhood $K$ of $y_0$ such that there exist $R > 1$ and $C_R > 0$ with the estimates:
$$ \vert\vert \omega_m\vert\vert_K \leq C_R R^m m! \quad \forall m \in \mathbb{N}.$$

\subsection{The $\Abc$-structure}

\begin{defn}\label{28/9 complex} 
Let $f : M \to \mathbb{C}$ be a holomorphic function and  denote the zero set of $f$  by $Y$. We assume that $f = 0$ is a reduced equation for $Y$. We define on $Y$ the following complexes of sheaves of left  $\Abc$-modules 
 $(\mathcal{K}_f^\bullet, D^\bullet)$ and $(\mathcal{I}_f^\bullet, D^\bullet)$ which are zero in negative degres as follows:
\begin{enumerate}
\item $\mathcal{K}_f^\bullet $ is the kernel of the map $\sum_{m=0}^\infty b^m\otimes \omega_m^\bullet \mapsto df\wedge \omega_0^\bullet $. Note that this kernel is an $\Abc$-submodule of the (graduate) $\Abc$-module 
$$ \oplus_{p=0}^{dim M} \  \Omega_{M\vert Y}^p\{\{b\}\} :=   \oplus_{p=0}^{dim M} \Big( \mathcal{O}_{M\vert Y}\{\{b\}\} \otimes_{\mathcal{O}_{M\vert Y}} \Omega_{M\vert Y}^p\Big) $$
where $ x \in \Abc$ acts on $g\otimes \omega \in \Big(\mathcal{O}_{M\vert Y}\{\{b\}\} \otimes_{\mathcal{O}_{M\vert Y}} (\Omega_{M\vert Y}^p)\Big)$ by the formula
$x[g \otimes \omega] = x[g]\otimes \omega$.
\item $\mathcal{I}_f^p$ is the sub-sheaf of the sheaf $\mathcal{K}_f^p$ defined by he condition $\omega_0 \in df\wedge \Omega_{M\vert Y}^{p-1}$ for each $p \geq 1$. 

\item  The differential $D^\bullet$ is defined by the formula
\begin{equation*}
 D^\bullet(\sum_{m=0}^\infty  b^m\omega_m^\bullet) = \sum_{m=0}^\infty b^m(d\omega_m^\bullet - df\wedge \omega_{m+1}^\bullet) \tag{Diff.}
 \end{equation*}
 \end{enumerate}
 \end{defn}
 
 Note that this definition contained several assertions which desserve to be proved !
 
 \parag{Proof of the assertions} To prove first that $D^p$ sends a germ at $y_0 \in Y$  of section of $\mathcal{K}_f^p$ to a germ at $y_0$ of section of $\mathcal{K}_f^{p+1}$
 we must show that the series in the right  hand-side of  Formula (Diff.) converges to a germ in  $\mathcal{K}_f^{p+1}$. This is consequence of the following facts:
 \begin{itemize}
 \item $d(Ker\, df^p) \subset Ker\, df^{p+1} \quad {\rm and} \quad df\wedge \Omega^p_{M\vert Y} \subset Ker\, df^{p+1} \quad   \forall p \geq 0.$
 \item For each compact neighborhood  $K'$  of $y_0$ contained in the interior of $K$ there exists a constant $C(K', K)$ such for any $p \geq 0$ and any  $\alpha \in H(K, \Omega_M^p)$ we have
 $$ \vert\vert df\wedge \alpha \vert\vert_{K'} \leq C(K', K)\vert\vert \alpha \vert\vert_K $$
 and 
 $$  \vert\vert d\alpha \vert\vert_{K'} \leq C(K', K) \vert\vert \alpha \vert\vert_K $$
 This simply means that the maps $\square \mapsto df \wedge \square$ and $\square \mapsto d\square$ from $H(K, \Omega_M^p)$ to $H(K', \Omega_M^{p+1})$ are ( linear and) continuous.
\end{itemize}
Then to see that $D^\bullet$ is $\Abc$-linear it is enough to remark that $D^\bullet$ commutes with $a$ and $b$. This is an easy consequence of the formulas:
  \begin{align*}
     & a\sum_{j=0}^{+\infty} b^j\omega_j  = f\omega_0 + \sum_{j=1}^{+\infty} b^j(f\omega_j + (j-1)\omega_{j-1}) \\
     & b\sum_{j=0}^{+\infty} b^j\omega_j  = \sum_{j=1}^{+\infty} b^j\omega_{j -1}
    \end{align*}

Finally the equalities $ab - ba = b^2$ and  $D^{p+1}\circ D^p = 0$ are easy computations left to the reader.$\hfill\blacksquare$\\

   We have a natural inclusion of complexes of left \ $\Abc$-modules
     $$\tilde{i} : (\mathcal{I}_f^\bullet, D^\bullet)  \to  (\mathcal{K}_f^\bullet, D^{\bullet})$$
     and a natural quotient map of complexes:
     $$ \tilde{j} : (\mathcal{K}_f^\bullet, D^{\bullet}) \to (\tilde{\mathcal{K}}_f^\bullet, D^\bullet).$$
     
     Remark that we have also  natural morphisms of complexes
       \begin{align*}
     & u :   (I^{\bullet}, d^{\bullet}) \to  (\mathcal{I}_f^\bullet, D^\bullet)  \\
     & v :  (K^{\bullet}, d^{\bullet}) \to  (\mathcal{K}_f^\bullet, D^{\bullet})\\
     & w:  (\tilde{K}^\bullet, d^\bullet) \to (\tilde{\mathcal{K}}_f^\bullet, D^\bullet)
    \end{align*}
    and that these morphisms are compatible with \ $i$,  $\tilde{i}$ and $j, \tilde{j}$. More precisely, this means that we have the commutative diagram of complexes
 $$   \xymatrix{0 \ar[d] & 0  \ar[d] \\ (I^{\bullet}, d^{\bullet}) \ar[d]^i \ar[r]^u &   (\mathcal{I}_f^\bullet, D^{\bullet}) \ar[d]^{\tilde{i}} \\
     (K^{\bullet}, d^{\bullet})  \ar[r]^v \ar[d]^j &  (\mathcal{K}_f^\bullet,  D^{\bullet})  \ar[d]^{\tilde{j}} \\
     (\tilde{K}^\bullet, d^\bullet)  \ar[r]^w \ar[d] & (\tilde{\mathcal{K}}_f^\bullet, D^\bullet) \ar[d]\\
     0 & 0 } $$
     
     where the vertical lines are exact sequence of complexes. So the fact that $u$ and $v$ are quasi-isomorphisms implies that $w$ is also a quasi-isomorphism.\\
     The proof that $u$ and $v$ are quasi-isomorphisms is the  goal of the next sub-section.

\subsection{The quasi-isomorphism}

 The following theorem is the convergent version  of Theorem 2.2.1. of \cite{[B.06]}.

\begin{thm}\label{29/9}
Let $M$ be a complex manifold and $f : M \to \mathbb{C}$ be a holomorphic function.  Denote the zero set of $f$  by $Y$ and assume that $f = 0$ is a reduced equation for $Y$.
The  natural  maps $u$, $v$  and $w$ 
\begin{align*}
& u : (I^\bullet, d^\bullet) \to (\mathcal{I}_f^\bullet, D^\bullet)  \\
& v: (K^\bullet, d^\bullet) \to (\mathcal{K}_f, D^\bullet) \\
& w:  (\tilde{K}^\bullet, d^\bullet) \to (\tilde{\mathcal{K}}_f^\bullet, D^\bullet)
\end{align*}
are quasi-isomorphims  of complexes,  compatible with the natural actions of $a$ and $b$ on the respective cohomology sheaves on $Y$. They are also compatible with the natural maps
$$ i :  (I^\bullet, d^\bullet) \to (K^\bullet, d^\bullet)  \quad {\rm and} \quad \tilde{i} : (\mathcal{I}_f^\bullet, D^\bullet) \to (\mathcal{K}_f, D^\bullet)  $$
$$ j : (K^\bullet, d^\bullet) \to (\tilde{K}^\bullet), d^\bullet) \quad {\rm and} \quad \tilde{j} :  (\mathcal{K}_f, D^\bullet) \to (\tilde{\mathcal{K}}_f^\bullet, D^\bullet). $$
\end{thm}

\parag{Proof} Denote by $(\widehat{K}_f^\bullet, D^\bullet)$ the formal completion\footnote{It should be noted that in a formal power series $\sum_{m=0}^\infty b^m\omega_m \in (\widehat{K}_f^\bullet)_{y_0}$ there exists, by definition (see footnote 3 in \cite{[B.08]}), an open neighborhood $U$ of  $y_0 \in M$ such that any  $\omega_m$, for each $m \geq 0$, belongs to  formal $f$-completion of $\Omega_M^\bullet(U)$.}  in $a$ and $b$ of the complex $(\mathcal{K}_f^\bullet, D^\bullet)$. We mean that we replace the sheaf $\mathcal{O}_{M\vert Y}\{\{b\}\}$ by $\mathcal{O}_{M\vert Y}[[b]]$ and the sheaf $(\Omega_{M\vert Y}^\bullet)$ by its formal completion in $f$. We have natural maps
$$ (K^\bullet, d^\bullet) \overset{\alpha}{\to}  (\mathcal{K}_f^\bullet, D^\bullet) \overset{\beta}{\to} (\widehat{K}_f^\bullet, D^\bullet) $$
where $\alpha$ is defined by $\alpha(u) = u$ and $D^\bullet(\alpha(u) = \alpha(du)$  for $u \in Ker\, df^\bullet$ and $\beta$ is the completion map. It is proved in \cite{[B.09]} Theorem 2.1.1 that $\beta\circ\alpha$ is a quasi-isomorphism. So $\alpha$ induces an injective map between cohomology sheaves. It is then enough to show that $\beta$ also induces an injective map between cohomology sheaves  to conclude that $\alpha$ is a quasi-isomorphism  (and then $\beta$ is also a quasi-isomorphism).\\
To prove that $\beta$ is  injective on the cohomology sheaves  consider $\omega$ a germ of section at the point $y_0 \in Y$ of the sheaf $\mathcal{K}_f^{p+1}$ which satisfies $D^{p+1}(\omega) = 0$ and  $\beta(\omega) = D^{p}(U)$ where $U$ is a section in $\widehat{\mathcal{K}}_f^{p}$ (so  satisfies $df \wedge u_0 = 0$). We shall represent $\omega$ by a series
$$ \omega = \sum_{m=0}^\infty b^m\omega_m $$
where $\omega_m$ are in $H(K, \Omega^{p+1}_M)$ with $K$ a compact neighborhood of $y_0$, satisfying the estimates $\vert\vert \omega_m\vert\vert_K \leq C_RR^m m!,\quad \forall m \in \mathbb{N}$. Put $U := \sum_{m=0}^\infty  b^mu_m $ where $u_m$ are in $\Omega^{p}_{M, y_0}$ with the condition $df \wedge u_0 = 0$. The relation $\beta(\omega) = D^{p}(U)$ gives the equalities of germs 
\begin{equation*}
 du_m - df\wedge u_{m+1} = \omega_m \quad \forall m \in \mathbb{N}.\tag{R}
 \end{equation*}
Define $V := \sum_{m=0}^\infty \frac{t^m}{m!} u_m $ in $\Omega^p_{M, y_0}[[t]]$ and $\Omega := \sum_{m=0}^\infty \frac{t^m}{m! }\omega_m$ in $H(L, \Omega^p_{M \times \mathbb{C}})$ where $L$ is the compact set in $M \times \mathbb{C}$ product of $K$ by $\{\vert t \vert \leq 1/2R\}$. The convergence of this series is insured by the estimates for the $\omega_m$. Now the differential of $V$ is given by
\begin{align*}
& dV = \sum_{m =1}^\infty \frac{t^{m-1}}{(m-1)!}dt\wedge u_m + \sum_{m=0}^\infty  \frac{t^m}{m!} du_m  \quad {\rm and \ so}\\
& dV = \sum_{m=0}^\infty  \frac{t^m}{m!} \big[dt\wedge u_{m+1} + df\wedge u_{m+1} + \omega_m \big]
\end{align*}
since $df\wedge u_0 = 0$ and the relations (R).
 This gives $dV - d(t + f)\wedge \partial_t V =   \Omega$.\\
 
 If we take local coordinates $x_1, \dots, x_{n+1}$ near the point $y_0$ in $M$, then near the point $(y_0, 0)$ in $M \times \mathbb{C}$ we may take local coordinates $x_1, \dots, x_{n+1}, t$ near $(y_0, 0)$ but also local coordinates $x_1, \dots x_{n+1}, z := t+f(x)$ near this point. Remark that $f(y_0) = 0$ since $y_0$ is in $Y$ and the point $(y_0, t=0)$  is equal to  the point $(y_0, z=0)$ in these new coordinates. Moreover, the vector field $\partial_t$ in the first coordinate system  is equal to the vector fiel $\partial_z$ in the second coordinate system. This allows to rewrite the formal equality above as
 $$ d_{/z}V = \Omega $$
 where $d_{/z}$ is the relative $z$-differential near $(y_0, 0)$, since we have, thanks to the equality $\partial_z = \partial_t$:
  $$   d_{/z} = d - dz\wedge \partial_z = d -  d(t+f)\wedge\partial_t $$ 
  for any  holomorphic (or formal)  germ of differential form at $(y_0, 0)$.\\
  So we have proved that the holomorphic germ $\Omega$ is  formally $d_{/z}$-exact at $(y_0, 0)$. But since $z$ is a smooth function (a coordinate !)
 the holomorphic germ $\Omega$ is $d_{/z}$-exact as a holomorphic germ. So we may find a convergent power series on a compact subset $K'\times \{\vert t\vert \leq \varepsilon\}$
 where $K' \subset\subset K $ is a compact neighborhood of $y_0$ and where $\varepsilon > 0$ is small enough such that
 $$ \tilde{V} := \sum_{m=0}^\infty  \frac{t^m}{m!} \tilde{v}_m $$
 satisfies  the estimates $\vert\vert \tilde{v}_m\vert\vert_{K'}\leq C_SS^m m! $ for $S > R$ large enough, and verifies on $K'\times \{\vert t\vert \leq \varepsilon\}$
 $$ d_{/z}(\tilde{V}) = \Omega \quad {\rm and \ so} \quad  d\tilde{V} = d(t+f)\wedge \partial_t(\tilde{V}) + \Omega.$$
 Then $\tilde{U} := \sum_{m=0}^\infty b^m\tilde{v}_m $ is a germ at $y_0$ of section of the sheaf $\mathcal{K}_f^{p}$ when $df\wedge \tilde{v}_0 = 0$ thanks to  the estimates given on the $\tilde{v}_m$.  It satisfies the equality
 $$ D^p(\tilde{U}) = \omega \quad {\rm  which \ implies} \quad d\tilde{v}_m - df \wedge \tilde{v}_{m+1} = \omega_m \quad \forall m \geq 0.$$
 
 But one point is missing now:  we dont know that $\tilde{v}_0$ is in $Ker\, df^{p}$ so that $\tilde{U}$ is may  not be  in $\mathcal{K}^p_f$.\\
 Note that for $p = 0$, $\tilde{v}_0 - v_0$ is a constant and changing $\tilde{U}$ by this constant gives the injectivity of $\beta$ on the cohomology in  degree $1$.\\
 This missing point for $p \geq 1$ is solved  thanks to the following result.
  \parag{Claim} There exists, for each  $p \geq 1$, an element  $w_0 \in \Omega_{M, y_0}^{p-1}$ such that $\tilde{v}_0 - v_0 = d w_0$ as elements in $\Omega_{M, y_0}^{p}$.
 
 \parag{Proof of the claim} Since $d_{/z}(\tilde{V} - V) = 0$ as a formal power series in $\Omega^{p}_{M, y_0}[[t]]$,  there exists $W \in \Omega^{p-1}_{M, y_0}[[t]]$ such that
 $$ \tilde{V} - V = d_{/z}W = dW - d(t+f)\wedge \partial_t W$$
 as a formal power series in $t$ since $z$ is a local coordinate. Looking at the terms without $t$ and $dt$  in the equality above
 we find that
 $$ \tilde{v}_0 - v_0 = dw_0 + df\wedge w_1\quad {\rm where} \quad W := \sum_{m=0}^\infty  \frac{t^m}{m!} w_m .$$
 This implies that $ df\wedge\tilde{v}_0 = df\wedge d w_0$. Then $\tilde{v}_0 - dw_0$ is in $K^p$ and $D^p(dw_0) = 0$. So $\tilde{U} - dw_0$ is in $\mathcal{K}_f^p$ and satisfies 
 $D^p(\tilde{U} - dw_0) = \omega$ proving that $\beta$ induces an injective map on the cohomology sheaves.\\

 The proof that $u$ is a quasi-isomorphism is obtained in the same way, replacing the condition $df\wedge \omega_0 = 0$ by $\omega_0 = df \wedge u'_0$ with $u'_0 \in \Omega^{p-1}$.\\
 Then the commutativity of the following  diagram 
 $$\xymatrix{0  \longrightarrow &\mathcal{E}_1\otimes \mathbb{C}_Y[1] \ar[d]^{\simeq}  \longrightarrow &(K^\bullet, d^\bullet) \ar[d]^j  \longrightarrow &(\tilde{K}^\bullet, d^\bullet) \ar[d]^{\tilde{j}}  \longrightarrow &  0 \\
 0  \longrightarrow &\mathcal{E}_1\otimes \mathbb{C}_Y[1]  \longrightarrow &(\mathcal{K}_f^\bullet, D^\bullet)  \longrightarrow &(\tilde{\mathcal{K}}_f^\bullet, D^\bullet)  \longrightarrow&  0} $$
 where the lines are exact sequences of complexes,  implies that $\tilde{j}$ is a quasi-isomorphism  since $j$ is a quasi-isomorphism.$\hfill \blacksquare$\\
 
 The main interest of the previous theorem is given in the following obvious Corollary:
 
 \begin{cor}\label{3/10} In the situation of Theorem \ref{29/9} the cohomology sheaves of the complexes $(I^\bullet, d^\bullet), (K^\bullet, d^\bullet), (\tilde{K}^\bullet, d^\bullet)$ have natural structure of left  $\Abc$-modules, deduced via the quasi isomorphisms $u, v, w$. Moreover, the natural maps $i$ and $j$ are left $\Abc$-linear.$\hfill\blacksquare$\\
 \end{cor}
 
 Remark that $ \tilde{K}$ has its cohomology sheaves supported in $S := \{df = 0 \}$.
 
\subsection{Generalized Brieskorn Modules associated to a holomorphic function}

The present sub-section is devoted to extend the results of \cite{[B.06]} Section 2 to the  sheaves of  left $\Abc$-modules constructed in the previous sub-section  on the relative De Rham cohomology   of a holomorphic function on a complex manifold.\\

We keep the situation of $f: M \to \mathbb{C}$ of  Theorem  \ref{29/9}.
  
  \begin{defn}\label{(a,b)-module}
  A {\bf generalized Brieskorn module (a GBM for short)} is  a left \ $\Abc$-module which is free and of finite rank on the commutative sub-algebra $B := \mathbb{C}\{\{b \}\}$ of $\Abc$.\\
  Since a GBM is a convergent (a,b)-module ($B[a]$ is a sub-algebra of $\Abc$; see \cite{[part I]}) we may extend standard properties of convergent (a,b)-modules to GBM.
So we say that a  generalized Brieskorn module $\mathcal{E}$ is :
  \begin{enumerate}
  \item  a {\bf  local} GBM  when \ $\exists N \in \mathbb{N}$ \ such that \ $a^N\mathcal{E} \subset b\mathcal{E}$;
  \item a {\bf simple pole} GBM  when \ $a\mathcal{E} \subset b\mathcal{E}$;
  \item a {\bf regular} GBM  when it is contained in a simple pole GBM.
    \item a {\bf geometric} GBM  when it is regular and  its Bernstein polynomial\footnote{Recall that the Bernstein polynomial of a regular GBM $\mathcal{E}$ is the minimal polynomial of the action of \ $-b^{-1}a$ \ on \ $\mathcal{E}^{\sharp}\big/ b\mathcal{E}^{\sharp}$ where $ \mathcal{E}^{\sharp}$ is  its saturation by $b^{-1}a$. So it is the Bernstein polynomial of $\mathcal{E}$ view as a convergent regular (a,b)-module. See \cite{[part I]}.} \ has its roots in \ $\mathbb{Q}^{-*}$ (this condition is motivated by the fundamental result of \cite{[K.76]}).
  \end{enumerate}
  \end{defn}

  Now let \ $\mathcal{E}$ \ be any left \ $\Abc$-module, and define \ $B(\mathcal{E})$ \ as the \ $b-$torsion of \ $\mathcal{E}$,  that is to say
  $$ B(\mathcal{E}) : = \{ x \in \mathcal{E} \ / \  \exists N \quad  b^Nx = 0 \}.$$
  Define \ $A(\mathcal{E})$ \ as the \ $a-$torsion of \ $\mathcal{E}$ \ and 
   $$\tilde{A}(\mathcal{E}) : = \{x \in \mathcal{E} \ / \ Bx \subset A(\mathcal{E}) \}.$$
   Remark that \ $B(\mathcal{E})$ \ and \ $\tilde{A}(\mathcal{E})$ \ are  left $\Abc$-submodules of \ $\mathcal{E}$ \ but that \ $A(\mathcal{E})$ \ is not stable by \ $b$ in general. 
  
  \begin{defn}\label{petit}
  A left \ $\Abc$-module \ $\mathcal{E}$ \ is {\bf small} when the following conditions hold
  \begin{enumerate}
 \item \ $B(\mathcal{E}) \subset \tilde{A}(\mathcal{E})$.
  \item \ $\exists N \ / \  a^N\tilde{A}(\mathcal{E}) = 0 $.
    \item \ $\mathcal{E}/B(\mathcal{E})$ \ is a finite type $B$-module.
  \end{enumerate}
  \end{defn}
  
 Moreover, in a  small $\Abc$-module  $\mathcal{E}$  we have always the equality \ $ B(\mathcal{E}) = \hat{A}(\mathcal{E})$.\\

  \begin{lemma}\label{crit. small}
   In a left  $\Abc$-module  $\mathcal{E}$ satisfying $1.$ and $2.$ of the previous definition,  we have always the equality \ $ B(\mathcal{E}) = \tilde{A}(\mathcal{E}) \subset Ker\, b^{2N}$.\\
  \end{lemma}

  \parag{Proof} 
    The condition $1.$ implies that the action of \ $b$ \ on \ $\tilde{A}(\mathcal{E})\big/B(\mathcal{E})$ \ is injective. But the condition $2.$ implies that \ $b^{2N} = 0$ \ on \ $\tilde{A}(\mathcal{E})$: \\
  From the commutation relation $ ab - ba = b^2 $ we deduce for each integer $N$  the identity  (see \cite{[B.06]}):
  \begin{equation*}
   N! b^{2N} = \sum_{j=0}^N  (-1)^j C^j_N b^j a^N b^{N-j}. \tag{F}
   \end{equation*}
   So we conclude that \ $\tilde{A}(\mathcal{E}) = B(\mathcal{E}) \subset Ker\, b^{2N}$.$\hfill \blacksquare$\\

 When $\mathcal{E}$ is a small $\Abc$-module,  the quotient \ $\mathcal{E}/B(\mathcal{E})$ \ is a GBM.

  \begin{defn}\label{geometric}
 When a   left \ $\Abc$-module \ $\mathcal{E}$ \ is  small,  the  GBM   \ $\mathcal{E}/B(\mathcal{E})$ \  is called the {\bf associated GBM} to $\mathcal{E}$.
   \end{defn}

   \begin{thm}\label{5/10} 
   The germ at each point $y_0 \in Y$ of the   cohomology sheaves of the complexes $ (K^\bullet, d^\bullet), (\tilde{K}^\bullet, d^\bullet)$ are {\bf small} left  $\Abc$-modules. The {\bf associated GBM  are geometric} and  in degree $p+1$ they  correspond respectively the { \bf nearby and vanishing cycle  systems} of the Gauss-Manin connection in degree $p$ of the germ $f_{y_0}: (M, y_0) \to \mathbb{C}, 0)$.
   \end{thm}

   \parag{Proof} We have to introduce some notations. Fix a point $y_0$ in $Y$ and denote by $\mathcal{H}^{p+1}_{y_0}$ the germ at $y_0$ of the  $(p+1)$-th cohomology
   sheaf of the complex $(K^\bullet, d^{\bullet})$.  Let $f_{y_0} : X_{y_0} \to D_{y_0}$ be a Milnor representative of the germ of $f$ at $y_0$ and denote by $F_{y_0}$ the Milnor fiber of $f$ at $y_0$. We fix a basis $\gamma_1, \dots, \gamma_d$ of $H^p(F_{y_0}, \mathbb{C})$ and denote by $\gamma_1^*, \dots, \gamma_d^*$ the dual basis in $H_p(F_{y_0}, \mathbb{C})$.\\
    Denote by $\mathscr{A} \in \mathbb{Q}/\mathbb{Z} \simeq ]0, 1] \cap \mathbb{Q}$  the image in $]0, 1] \cap \mathbb{Q}$ of the opposite of the roots of the Bernstein polynomial of  the germ $f_{y_0}$. We shall consider the  GBM  (see  the example following Theorem \ref{div.3})
   $$ \Xi^{p}_{\mathscr{A}} \otimes H^p(F_{y_0}, \mathbb{C})$$ 
   which is a simple pole geometric  GBM containing the integrals 
   $$ \int_{\gamma_s} \omega/df $$
   when $\omega$ is in $K^{p+1}_{y_0} \cap Ker \, d$ and $\gamma_s$ is a horizontal family of compact  $p$-cycles in the fibers of $f_{y_0}$.\\
   \parag{Claim}There is   a $\Abc$-linear map
   $$ \Phi : \mathcal{H}^{p+1}_{y_0}  \to \Xi^{p}_{\mathscr{A}} \otimes H^p(F_{y_0}, \mathbb{C})$$ 
   given by $\omega \mapsto \sum_{j=1}^d \big(\int_{\gamma^*_{j, s}} \omega/df \big) \otimes \gamma_j $ whose kernel is the $b$-torsion of $\mathcal{H}^{p+1}_{y_0}$.
   
   \parag{Proof of the Claim} This map is well defined and its $\Abc$-linearity is an easy consequence of the relations $\Phi(a\omega) = a\Phi(\omega)$ and $\Phi(b\omega) = b\Phi(\omega)$ which are clear.\\
   Since the target has no $b$-torsion, the kernel contains the $b$-torsion of $\mathcal{E} := \mathcal{H}^{p+1}_{y_0}$. Conversely, if $\omega$ is in the kernel of $\Phi$ then,  the $d$-closed $p$-form $\omega/df$ on $F_{y_0}$ induces $0$ in $H^p(F_{y_0}, \mathbb{C})$. This easy to see using a Jordan basis for the action of the monodromy of $f_{y_0}$ on $H^p(F_{y_0}, \mathbb{C})$. So, thanks to \cite{[G.65]} there exists a germ at $ y_0$ of a meromorphic $(p-1)$-form $u$ with poles in $Y$ such that $\omega = df\wedge du$. Now writing 
   $u = f^{-N}v$ where $v$ is in $\Omega_{y_0}^{(p-1)}$ we obtain
   $$ f^N\omega = df \wedge f^Nu = df\wedge dv$$
   which means that $a^N[\omega] = 0$ in $\mathcal{H}^{p+1}_{y_0}$.\\
   So the kernel is also contained in the $a$-torsion, so in $\tilde{A}(\mathcal{E})$. \\
   But since $\tilde{A}(\mathcal{E})$ is clearly contained in $Ker\, \Phi$ we obtain $B(\mathcal{E})= Ker\, \Phi = \tilde{A}(\mathcal{E})$.\\
   Since for $y_0$ given, there exists $\kappa \in \mathbb{N}$ such that $f^\kappa (Ker \, df^{\bullet+1})_{y_0} \subset df\wedge \Omega^\bullet_{y_0}$ we have
   $a^\kappa \mathcal{E} \subset b\mathcal{E}$ and then, thanks to Formula $(F)$,  $b^{2\kappa}B(\mathcal{E}) = 0$. Then we conclude that for $N = 4\kappa^2$ we have
    $a^N\tilde{A}(\mathcal{E}) = 0$ proving the property 2. of a small $\Abc$-module for $\mathcal{E}$.\\
      Now   remark that, since $\Xi^{p}_{\mathscr{A}} \otimes H^p(F_{y_0}, \mathbb{C})$ is a free finite type $B$-module, its $\Abc$- sub-module isomorphic to
   $\mathcal{H}^{p+1}_{y_0} \big/(b-torsion)$ is also a free finite type $B$-module so it is  a geometric  GBM.\\
   The case of the complex $(\tilde{K}^\bullet, d^\bullet)$ is analogous using only the singular part of the asymptotic expansions, so killing $(\mathbb{C}\{s\} \otimes V) $ in $(\Xi_\mathscr{A}^n \otimes V)$.$\hfill \blacksquare$
   
   \parag{Remark} The exact sequence  of complexes\footnote{ $\mathcal{E}_1$ is placed in degree 1 and  $\{0\}$ is in all other degrees.} 
   $$0 \to (\mathcal{E}_1)_{y_0} \to (K^\bullet, d^\bullet)_{y_0} \to (\tilde{K}^\bullet, d^\bullet)_{y_0} \to 0 $$
   corresponds to the canonical map between nearby and vanishing cycles.

   \subsection{The  global  finiteness Theorem}
  
  The main result of this sub-section is the following theorem, which shows that the global  Gauss-Manin connection of a proper holomorphic function produces geometric GBM associated to vanishing cycles and nearby cycles.
  In fact the compactness of the singular fiber is enough to prove the finiteness result for the GBM associated to the global vanishing cycles in suitable situations.
  
  Introduce now the complex $(\tilde{I}^\bullet, d^\bullet)$ via the exact sequence of complexes on $Y_0$:
  
  $$ 0 \to \mathcal{E}_1\otimes \mathbb{C}_Y \to (I^\bullet, d^\bullet) \to (\tilde{I}^\bullet, d^\bullet) \to 0 $$
  
  where $\mathcal{E}_1\otimes \mathbb{C}_Y$ is in degree $1$ and is the only non zero term in this complex.
  It is clear that for $p \geq 1$ the natural map  $\mathcal{H}^{p+1}(I^\bullet, d^\bullet) \to \mathcal{H}^{p+1}(\tilde{I}^\bullet, d^\bullet) $ is an isomorphism and that in degre $1$
  we have the exact sequence 
  $$0 \to  \mathcal{E}_1\otimes \mathbb{C}_Y \to  \mathcal{H}^{1}(I^\bullet, d^\bullet) \to  \mathcal{H}^{1}(\tilde{I}^\bullet, d^\bullet) \to 0.$$
  
  Then we have an exact sequence of complexes
  \begin{equation*}
  0 \to (\tilde{I}^\bullet, d^\bullet) \to (\tilde{K}^\bullet, d^\bullet) \to (L^\bullet, d^\bullet) \to 0 \tag{@}
  \end{equation*}
  where $L^p = Ker \, df^p\big/Im \, df^p $ is a coherent sheaf on $S$ the singular set of $Y$. \\
  
  \begin{thm}\label{finite} Let $f : M \to \mathbb{C}$ be a non constant holomorphic function on a complex manifold $M$. Assume that $Y := \{f = 0 \}$ is reduced  and contains the set $S := \{df = 0\}$  of $f$. Our hypohesis is now that $S$ is compact. Then the hyper-cohomology $\mathbb{H}^{p+1}(Y, (\tilde{K}^\bullet, d^\bullet))$ is a left   $\Abc$-module for each $p \geq 0$  such that its quotient by its $b$-torsion is a free finite type $B$-module. Moreover its $b$-torsion is finite dimensional.
  \end{thm}

  \parag{Proof} The first point is to remark that for each $p \geq 0$ we have an isomorphism of left the $\Abc$-modules $ \mathcal{H}^{p+1}(\tilde{I}^\bullet, d^\bullet) \simeq b\mathcal{H}^{p+1}(\tilde{K}^\bullet, d^\bullet) $:\\
 If  $\omega \in \Omega^{p+1}_{y_0}$ satisfies $df \wedge\omega = 0$ and $d\omega = 0$ then $b[\omega] $ in $\mathcal{H}^{p+1}(\tilde{K}^\bullet, d^\bullet)$ is induced by $df \wedge u$ \  if \ $du = \omega$. So it belongs to the image of $ \mathcal{H}^{p+1}(\tilde{I}^\bullet, d^\bullet)_{y_0}$ in $\mathcal{H}^{p+1}(\tilde{K}^\bullet, d^\bullet)_{y_0}$.\\
  Conversely if $df \wedge u$ is $d$-closed then its image in $\mathcal{H}^{p+1}(\tilde{K}^\bullet, d^\bullet)$ is $b[du]$ where $du$ is in $Ker\, df^{(p+1)}$. This proves our  first assertion. \\
 Let us  show that \ $\mathcal{E} : =  \mathbb{H}^{p+1}(Y, (\tilde{K}^{\bullet}, d^{\bullet}))$ \ satisfies  conditions in the statement.  \\
   Consider now the long exact sequence of hyper-cohomology of the exact sequence of complexes
 $$   0 \to (\tilde{I}^{\bullet}, d^{\bullet}) \to (\tilde{K}^{\bullet}, d^{\bullet}) \to ([K/I]^{\bullet}, d^{\bullet}) \to 0 .$$
 It contains the exact sequence
 $$  \mathbb{H}^{p-1}(Y, ([K\big/I]^{\bullet}, d^{\bullet})) \to \mathbb{H}^p(Y, (\tilde{I}^{\bullet}, d^{\bullet})) \overset{\mathbb{H}^p(i)}{\to} \mathbb{H}^p(Y, (\tilde{K}^{\bullet}, d^{\bullet})) \to \mathbb{H}^{p}(Y, ([\hat{K}\big/\hat{I}]^{\bullet}, d^{\bullet})) $$
 and we know that \ $b$ \ is induced on the complex of \ $\Abc$-modules  quasi-isomorphic to \ $(\tilde{K}^{\bullet}, d^{\bullet})$ \ and  $(\tilde{I}^{\bullet}, d^{\bullet})$  by the composition \ $i\circ \tilde{b}$ \ where \ $\tilde{b}$ \ is a quasi-isomorphism of complexes of \ $B$-modules. This implies that the kernel  and the co-kernel of \ $\mathbb{H}^p(i)$ \ are isomorphic (as \ $\mathbb{C}$-vector spaces) to \ $Ker\, b$ \ and \ $Coker \, b$ \ respectively. So  it is enough to prove finite dimensionality for the vector spaces \ $  \mathbb{H}^{j}(Y, ([K\big/I]^{\bullet}, d^{\bullet})) $ \ for all \ $j \geq 0 $. But the sheaves \ $[K\big/I]^j \simeq [Ker\,df\big/Im\, df]^j$ \ are coherent on \ $M$ \ and supported in \ $S$. The spectral sequence
 $$ E_2^{p,q} : = H^q\big( H^p(M, [K\big/I]^{\bullet}), d^{\bullet}\big) $$ 
 which converges to \ $ \mathbb{H}^{j}(M, ([K\big/I]^{\bullet}, d^{\bullet})) $, is a bounded complex of finite dimensional vector spaces by Cartan-Serre Theorem since $S$ is compact.
  This gives the desired finite dimensionalities. $\hfill \blacksquare$\\
  
  The following Corollary is an  example of consequence of the previous Theorem when we are in a situation where the local system $Rf_*\mathbb{Z} $ has a quasi-unipotent monodromy and defines a regular Gauss-Manin connection at the origin. The proof is similar to the proof of Theorem \ref{5/10} by considering the quasi-isomorphic complex $(K^\bullet_\infty , d^\bullet)$ of $ \mathscr{C}^\infty$ differential forms  $\omega$ such that $df \wedge \omega = 0$ which is quasi-isomorphic to $(K^\bullet, d^\bullet)$ thanks to Lemma 6.1.1 in \cite{[B.08]},  to obtain, since $K^\bullet_\infty$ is a fine sheaf
$$  H^{p+1}(H^0(Y, K^\bullet_\infty) , d^\bullet) \simeq   \mathbb{H}^{p+1}(Y, (K^{\bullet}, d^{\bullet})) $$
  
and the asymptotics of the  global  period-integrals given by integration of $d$-closed global sections  on $Y$ of the sheaf $K_\infty^\bullet$ to obtain an injective (after quotient by the $b$-torsion which contains the $a$-torsion by the same argument as in Theorem \ref{5/10}) $\Abc$-linear map to a geometric GMB of the type $\Xi_\mathscr{A}^N \otimes H^\bullet(Y_1, \mathbb{C})$ where $Y_1$ is a generic fiber of the map $f$.\\
  
  \begin{cor}\label{6/10} Let $f : M \to D$ be a proper holomorphic map between a connected complex manifold of dimension $n+1, n \geq 1$, such that the subset $\{df = 0 \}$ is contained in $Y := \{f = 0\}$. We assume that $\{f = 0\}$ is  a reduced equation of $Y$. Then the global hyper-cohomology on $Y$ of the complexes $(K^\bullet, d^\bullet)$ and $(\tilde{K}^\bullet, d^\bullet)$ are small left $\Abc$-modules such that their associated GBM are geometric. Moreover, for each $p \geq 1$ the canonical map
  $ can : \mathbb{H}^p(Y, (K^\bullet, d^\bullet)) \to \mathbb{H}^p(Y, (\tilde{K}^\bullet, d^\bullet) $
  is $\Abc$-linear.$\hfill \blacksquare$
  \end{cor}

  \parag{Some questions}  Let $f : M \to D$ a proper holomorphic map of a   Kahler (or quasi-projective)  manifold $M$ onto a disc such that the set  $\{df = 0\}$ is contained in $Y := \{f = 0\}$ which is assumed to be a  reduced hyper-surface.\\
 Precise the relations between the mixte Hodge structure on the cohomology of the fiber $\{f = 0\}$ and the generalized Brieskorn modules which are associated to the hyper-cohomology of the complex  $(Ker\,df^\bullet, d^\bullet)$ on $Y$.\\
 More precisely, denoting $H^p(Y)$ the $p$-th cohomology group of $Y$ and $\mathcal{H}^{p+1}$ the  generalized Brieskorn module associated to the $(p+1)$-th hyper-cohomology group of the complex  $(Ker\,df^\bullet, d^\bullet)$ on $Y$, precise the relations:
  \begin{enumerate}
 \item  between the Hodge filtration on $H^p(Y)$ and the $b$-filtrations on $\mathcal{H}^{p+1}$. (the natural one and the one induced by the $b$-filtration of the saturation of $\mathcal{H}^{p+1}$);
 \item between the weight filtration  on $H^p(Y, \mathbb{Q})$ and the semi-simple filtration on $\mathcal{H}^{p+1}$ introduced in \cite{[part I]};
 \item between the higher Bernstein polynomials of $\mathcal{H}^{p+1}$ introduced in \cite{[part II]} and  the SHM on $H^p(Y)$ (see \cite{[D.71]})  or the mixed Hodge module defined by Morihiko Saito \cite{[S.90]}.
 \end{enumerate}


\newpage

  
    \section*{Bibliography}
  
  \begin{thebibliography}{99}
 


\bibitem{[Br.70]} Brieskorn, E. {\it Die Monodromie der Isolierten Singularit{\"a}ten von Hyperfl{\"a}chen},\\
 Manuscripta Math. 2 (1970), pp. 103-161.

\bibitem{[B.81]} Barlet, D.  \textit{D\'eveloppements asymptotiques des fonctions obtenues par int\'egration sur les fibres}, 
 Inv. Math. vol. 68 (1982), pp. 129-174.

\bibitem{[B.84-a]} Barlet, D. \textit{Contribution effective de la monodromie aux  d\'eveloppements \\ asymptotiques}, 
 Ann. Ec. Norm. Sup. t.17 (1984), pp. 293-315.

\bibitem{[B.84-b]}  Barlet, D. \textit{Contribution du cup-produit de la fibre de Milnor aux p\^oles de \ $\vert f \vert^{2\lambda}$},  Ann. Inst. Fourier (Grenoble) t. 34, fasc. 4 (1984), pp. 75-107.

\bibitem{[B-M.89]}  Barlet, D. et Maire, H.M. \textit{Asymptotic expansion of complex integrals via Mellin transform}, 
 J.  Funct. Anal. 83 (1989), pp. 233-257.


\bibitem{[B.93]} Barlet, D. {\it Th\'eorie des (a,b)-modules I},\\
 in Complex Analysis and Geometry, Plenum Press, (1993), pp. 1-43.
 
 \bibitem{[B.97]} Barlet, D. {\it Th\'eorie de (a,b)-modules II. Extensions.}\\
 in Complex Analysis and Geometry, Pitman  Notes $n^0 366$, (1997) pp. 19-59.

\bibitem{[B.06]} Barlet, D. {\it Sur certaines singularit\'es non isol\'ees d'hypersurfaces I},\\
 Bull. Soc. math. France 134 (2), (2006), pp.173-200.

\bibitem{[B.08]} Barlet, D. {\it Sur certaines singularit\'es non isol\'ees d'hypersurfaces  II}, \\
J. Alg. Geom. 17 (2008), pp. 199-254.

\bibitem{[B.09]} Barlet, D. {\it P\'eriodes \'evanescentes et (a,b)-modules monog\`enes}, \\
Bollettino U.M.I. (9) II (2009), pp. 651-697.

\bibitem{[B.III]} Barlet, D. {\it Sur les fonctions \`a lieu singulier de dimension 1 }, \\
Bull. Soc. math. France 137 (4), (2009),  pp. 587-612.

 
 \bibitem{[B.22]} Barlet, D. {\it Algebraic differential equations of period-integrals},\\
   Journal of Singularities Volume 25 (2022), pp. 54-77.
 
 \bibitem{[B.23]} Barlet, D. {\it Complement to Higher Bernstein Polynomials and Multiple \\
 Poles of ...} arXiv:2311.13259 math. AG

\bibitem{[B.-S. 04]} Barlet, D. et Saito, M. {\it Brieskorn modules and Gauss-Manin systems for non isolated hypersurface singularities,} \\
J. Lond. Math. Soc. (2) 76 (2007) $n^01$ \ pp. 211-224.

\bibitem{[part I]} Barlet, D. {\it Generalized Bernstein Modules I; convergent (a,b)-modules},\\
 hal-04972623v2.

\bibitem{[part II]} Barlet, D. {\it Generalized Bernstein Modules II; Higher Bernstein polynomials and Multiple poles}, arXiv:2503.04383.

\bibitem{[Bj.93]} Bj{\"o}rk, J-E, {\it   Analytic D-modules and applications}, Kluwer Academic  publishers (1993).

\bibitem{[D.70]} Deligne, P. {\it Equations diff\'erentielles \`{a}  Points Singuliers R\'eguliers}\\
 L.N. in Maths 163 (1970) Springer Verlag.
 
 \bibitem{[D.71]} Deligne, P. {\it Th\'eorie de Hodge II}\\
 Publ. IHES 40 (1971) pp. 5-57.

\bibitem{[G.65]} Grothendieck, A. {\it On the de Rham cohomology of algebraic varieties}, \\
Publ. Math. IHES 29 (1966), pp. 93-101.

 \bibitem{[K.76]} Kashiwara, M. {\it b-function and holonomic systems},\\
  Inv. Math. 38 (1976) pp. 33-53.
 

\bibitem{[M.74]} Malgrange, B. {\it Int\'egrale asymptotique et monodromie}, \\
Ann. Sc. Ec. Norm. Sup. 7 (1974), pp. 405-430.

\bibitem{[M.75]} Malgrange, B. {\it Le polyn\^ome de Bernstein d'une singularit\'e isol\'ee},\\
 in Lect. Notes in Math. 459, Springer (1975), pp. 98-119.

\bibitem{[K.S]} Saito, K. {\it Period mapping associated to a primitive form}\\
Publ. RIMS 19  (1983) $n^03$, pp. 1231–1264. 


\bibitem{[S.89]} Saito, M. {\it On the structure of Brieskorn lattices}, \\
Ann. Inst. Fourier 39 (1989), pp. 27-72.

\bibitem{[S.90]}  Saito, M. {\it Mixtes  Hodge modules} ,\\
Publ. RIMS. 26 (1990), pp.221-333


\end{thebibliography}
\end{document}